\documentstyle{amsart}
\author{Rogier Brussee}
\date{alg-geom/9503004}

\title[$C^\infty$ properties of K\"ahler surfaces]
{The  canonical class and the  $C^\infty$-properties of K\"ahler surfaces.}

\address{Fakult\"at f\"ur Mathematik Universit\"at Bielefeld \\
Postfach 100131\\
33501 Bielefeld}
\email{brussee@@mathematik.uni-bielefeld.de}

 \makeatletter

\def\defmathcal{\def\mathcal##1{{\cal ##1}}}
\def\defmathrm{\def\mathrm##1{{\rm ##1}}}
\def\defmathit{\def\mathit##1{{\it ##1}}}
\def\defBbb{\def\Bbb##1{{\bf ##1}}}
\def\defmathbb{\def\mathbb##1{{\Bbb ##1}}}

\@ifundefined{mathcal}{\defmathcal}{}
\@ifundefined{mathrm}{\defmathrm}{}
\@ifundefined{mathit}{\defmathit}{}
\@ifundefined{Bbb}{\defBbb}{}
\@ifundefined{mathbb}{\defmathbb}{}

\def\oldlatexwarning{
\hbox to 0pt{
$\phantom{\mathcal WHEN THIS GIVES AN ERROR, UNCOMMENT FIRST
DEFINITION}
\phantom{
\mathrm  WHEN THIS GIVES AN ERROR, UNCOMMENT SECOND
DEFINITION}
\phantom{\mathit WHEN THIS GIVES AN ERROR, UNCOMMENT THIRD
DEFINITION}$
\hskip 0pt minus 1fil
}
}

 \theoremstyle{plain}

 \begingroup
 \theorembodyfont{\sl}

 \newtheorem{Theorem}{Theorem}
 \newtheorem{Corollary}[Theorem]{Corollary}
 \newtheorem{Proposition}[Theorem]{Proposition}
 \newtheorem{Lemma}[Theorem]{Lemma}
 \endgroup

 \theoremstyle{definition}

 \newtheorem{Definition}[Theorem]{Definition}
 \newtheorem{Properties}{Properties}

 \theoremstyle{remark}

 \newtheorem{Remark}[Theorem]{Remark}

 \newtheorem{Acknowledgment}{Acknowledgment}

 \def\txt#1{\hbox{ #1 }}

 \newdimen\arrowlen
 \def\longmap{\arrowlen3em }
 \def\m@p#1#2#3#4{
	\buildrel
	       \hbox spread \arrowlen{\skip@=-0.33\arrowlen plus 1 fil
		     \hskip\skip@$\m@th\scriptstyle #4$\hskip\skip@
	       }
	\over{
	       \mathord#1\mkern-6mu
	       \cleaders\hbox{$\mkern-2mu\mathord#2\mkern-2mu$}\hfill
	       \mkern-6mu\mathord#3
	}
 }

 \def\map{\m@p--\rightarrow}
 \def\pam{\m@p\leftarrow--}
 \def\equal{\m@p===}
 \def\Map{\m@p==\Rightarrow}
 \def\Pam{\m@p\Leftarrow==}
 \def\vm@p#1#2{\Big#1
	\rlap{$\vcenter{\hbox{$\scriptstyle #2$}}$}}
 \def\downmap{\vm@p\downarrow}
 \def\upmap{\vm@p\uparrow}
 \def\vequal{\vm@p\Vert}
 \def\Downmap{\vm@p\Downarrow}
 \def\maps{\mapstochar\nobreak}
 \def\inj{\lhook\nobreak\joinrel\nobreak}
 \def\surj#1{#1\nobreak\mathrel{{\mkern-12mu}}\nobreak\rightarrow}

\def\rmmath#1{\mathop{\mathrm{#1}}\nolimits}

\def\Hom{\rmmath{Hom}}

\def\End{\rmmath{End}}

\def\sheafHom{\mathop{{\mathcal H}\mkern-3mu{\mathit om}}\nolimits}

\def\sheafEnd{\mathop{{\mathcal E}\mkern-3mu{\mathit nd}}\nolimits}

\def\invlim{{\displaystyle \lim_{\leftarrow}}}

\def\Pic{\rmmath{Pic}}

\def\Proj{\rmmath{Proj}}

\def\rank{\rmmath{rank}}
\def\kod{\kappa}

\def\Ker{\rmmath{Ker}}
\def\Im{\rmmath{Im}}
\def\Coker{\rmmath{Coker}}

\let\tensor=\otimes
\def\directsum{\mathop\oplus}

\def\ch{\rmmath{ch}}
\def\td{\rmmath{td}}

\let\slant=/

\mathcode`\:="603A 
\mathchardef\:="303A 

\let\iso=\cong
\def\eqdef{\buildrel\mathrm{def}\over=}
\def\union{\cup}
\def\supsupset{\supset\!\supset}
\def\subsubset{\subset\!\subset}

\def\id{\rmmath{id}}

\def\R{{\Bbb R}}
\def\C{{\Bbb C}}
\def\Q{{\Bbb Q}}
\def\Z{{\Bbb Z}}
\def\N{{\Bbb N}}
\def\P{{\Bbb P}}

\def\O{{\mathcal O}}

\def\numfrac#1#2{\mathchoice{{\textstyle{ #1\over#2}}}%
{{ #1\over#2}}{{#1/#2}}{{#1/#2}}}
\def\quart{{\numfrac14}}
\def\half{{\numfrac12}}

\let\next=\~
\let\~=\tilde
\let\tilde=\next

\let\next=\^
\let\^=\hat
\let\hat=\next

\def\"#1{{\accent"7F \if#1i\i\else#1\fi}}
\def\o#1{\mathord{\buildrel\circ\over#1}}

\def\dual{{\scriptscriptstyle\vee}}

\def\({\left(}
\def\){\right)}

\def\U{\mathrm{U}}

\def\E{{\mathcal E}}
\def\L{{\mathcal L}}
\def\M{{\mathcal M}}
\def\<#1>{\left<#1\right>}

\newif\ifcomment
\def\comment#1{\ifcomment #1 \fi}

\newcommand\nc\newcommand

\nc\Spin{{\mathrm{Spin}}}
\nc\delbar{\rlap{/}\partial}
\nc\dbar{\bar\partial}
\nc\Kmin{K_{\min}}
\nc\Xmin{X_{\min}}
\nc\K{{\mathcal K}}
\nc\NE{ {\rmmath{NE}}}
\nc\NEbar{\overline{\NE}}

\nc\eps{\varepsilon}
\nc\sign{\sigma}
\nc\SO{\mathrm{SO}}
\nc\SC{{\mathcal{SC}}}
\nc\G{{\mathcal G}}
\nc\A{{\mathcal A}}
\nc\B{{\mathcal B}}
\nc\Pee{{\mathcal P}}
\nc\W{{\mathcal W}}
\nc\MM{{\Bbb M}}
\nc\MBN{\M^{\mathrm{BN}}}
\nc\cl{{\mathrm cl}}
\nc\Ind{{\rmmath{Ind}}}
\nc\eN{{\mathcal N}}
\nc\Fred{\mathrm{Fred}}
\nc\laplace{\mathop{\Delta}}
\nc\Vol{\rmmath{Vol}}
\nc\Pbar{{\bar{\P}}}
\nc\dee{\partial}
\nc\vdim{\rmmath{vdim}}
\def\Re{\rmmath{Re}}
\nc\Cee{{\mathcal C}}
\nc\Herm{{\mathcal H}}
\nc\Amod{\A^{\mathrm{mod}}}
\nc\PeeBN{{\Pee^{\mathrm{BN}}}}
\nc\locsys{\xi}
\nc\orr{or}
\nc\cH{\check H}
\nc\hol{{\mathrm{hol}}}
\nc\loc{{\mathrm{loc}}}
\nc\harm{{\mathrm{harm}}}
\nc\Top{{\mathrm{top}}}
\nc\sfA{\mathord{\hbox{\sf A}}}
\nc\ddt{{d \over dt}}
\makeatother

\begin{document}
\maketitle

\begin{abstract}
\oldlatexwarning
We give a self contained proof	that
for K\"ahler surfaces with non negative Kodaira dimension,
the canonical class of the minimal model and the
$(-1)$-curves, are oriented diffeomorphism invariants up to sign.
This includes the case $p_g = 0$.
It
implies that the Kodaira dimension is determined by the underlying
differentiable manifold.
We compute the Seiberg Witten invariants of all K\"ahler
surfaces of non negative Kodaira dimension.
We then reprove that the multiplicities of the elliptic
fibration are determined by the underlying oriented manifold, and that
the plurigenera of a surface are oriented diffeomorphism invariants.
The proof uses a set up of Seiberg Witten theory that replaces generic
metrics by the	construction of a localised Euler class of an infinite
dimensional bundle with a Fredholm section. This  makes the
techniques of excess intersection available in gauge theory
\end{abstract}

A compact complex surface $X$ with non negative Kodaira dimension, has a
unique minimal model $\Xmin$. The pullback of the canonical line bundle
minimal model $\omega_{\min}$ is in some ways the most basic birational
invariant of the surface, if only because it is the polarisation $\O(1)$
of the canonical model $\Proj(\directsum H^0(nK))$.  Recently,
Kronheimer, Mrowka and Tian,  Yau proved that the cohomology class,
$\Kmin = c_1(\omega_{\min})\in H^2(X,\Z)$ is invariant under oriented
diffeomorphism up to sign for minimal surfaces of general type
with $p_g >0$ \cite{Stern:talk}.  While completing this manuscript,
Friedman and Morgan posted a proof for the case $p_g =0$ \cite{FM:SW}.
In the case of elliptic surfaces it was already known to be true by the
joint effort of many people,  as it is a direct consequence of the
invariance of the multiplicities of the elliptic fibration.

The proof is  based on fundamental work of Witten and Seiberg
\cite{Witten}, who introduced a new set of non linear equations, the
monopole equations, which allow to define new differentiable invariants,
similar in spirit to the Donaldson invariants, but much easier to
handle. During the Stillwater conference in November 1994, Stefan Bauer
and I learned about these new invariants and the invariance of $\Kmin$,
and we decided to run a seminar on Seiberg Witten theory in Bielefeld.
I started to think about a proof for the invariance of $\Kmin$ using the
SW invariants along the lines of \cite{ikke:KM}. A little to my
surprise, it worked out beautifully. The monopole  equations which
define the  SW classes, once specialised to the K\"ahler case, give all
the necessary information, without the necessity to  prove a full Thom
conjecture type of result.

In fact it even turned out to be possible to deal with the case $p_g =
0$  and the elliptic surfaces simultaneously and almost uniformly.
There is no need to use any result from ``classical '' Donaldson theory
either. From the point of view of classification of surfaces,  it is
very satisfactory that the various levels of nefness of $\Kmin$  (nef
and big, nef but not big, torsion)  is what makes the proof work for
Kodaira dimension $\kappa \ge 0$, what makes it fail for the rational
and ruled case, and what makes for the difference in the different
Kodaira dimensions. If $p_g = 0$ the higher plurigenera, and in
particular the 2-canonical system plays a direct role.

While proving the invariance of the canonical class, we have to prove
the invariance of  $(-1)$-curves as well. This leads directly to the
differentiable characterisation corollary~\ref{-infty}	of rational and
ruled surfaces which are characterised algebraically by the existence of
a smooth rational curve $l$  with $l^2 \ge 0$ \cite[Prop. V.4.3]{BPV}.
The invariance of the Kodaira dimension (Van de Ven conjecture) and the
invariance of the plurigenera for surfaces of general type is then
almost an afterthought. The Van de Ven Conjecture had already been
solved	with  Donaldson theory (\cite{FM} for all surfaces but rational
surfaces and surfaces of general type with $p_g = 0$, and Friedman Qin
\cite{FQ} and Pidstrigatch \cite{PT},\cite{Viktor:easy blowup} for the
remaining cases, see also \cite{OkonekTeleman:VdVen} for an easy proof
of the remaining case with Seiberg Witten theory).

\begin{Theorem}\label{main}
If  $X$ is a K\"ahler  surface of non negative Kodaira
dimension then
\begin{enumerate}
\item \label{maina}
$\Kmin$ is determined by the underlying oriented manifold  up to sign,
\item \label{mainb}
every $(-1)$-sphere is $\Z$-homologous to  a $(-1)$-curve up to sign.
\end{enumerate}
\end{Theorem}

\comment{
Actually the statements are true on non-K\"ahler surfaces by
classification ($K$ is trivial for minimal primary Kodaira surfaces, and
for secondary Kodaira surfaces the class of fibre and the order of $K$
is determined by the fundamental group, moreover minimal
Kodaira surfaces are $K(\pi,1)$'s , so the $(-1)-curves are the image of
$\pi_2$)
}

\begin{Corollary}{}\label{spheres}
If a K\"ahler surface $X$ has non negative Kodaira dimension then
every smooth sphere $S$ with $S^2 \ge 0$ is $\Z$-homologous to $0$,
\end{Corollary}
\begin{Corollary}\label{-infty}
A  K\"ahler surface is rational or ruled if and only if it contains a
smooth sphere $S \ne 0 \in H^2(X,\Z)$ with $S^2 \ge 0$.
\end{Corollary}
\begin{Corollary}\label{Kodaira}
The Kodaira dimension of a K\"ahler surface
is determined by the underlying differentiable manifold.
\end{Corollary}

To deal with the case $p_g =0$, we encounter higher dimensional
moduli spaces and more to the point, moduli spaces that have larger
than virtual dimension. However, following Pidstrigatch and Tyurin,
we will identify the multiplicity of a Seiberg
Witten class  as a  localised Euler class of an infinite rank bundle
with a Fredholm section, and the oversized  moduli spaces will cause
no problem at all.
In addition without much extra work, the computation will give us
the Seiberg Witten multiplicities of all K\"ahler elliptic surfaces.
It demonstrates my belief that Seiberg Witten
theory for surfaces is completely computable.
An elegant argument of Stefan Bauer,
then gives yet another proof that for elliptic surfaces with finite cyclic
fundamental group the multiplicities of the elliptic fibration are
determined by the underlying oriented manifold. The oriented homotopy type
determines the multiplicities for other elliptic surfaces
\cite[Theorem S.7]{FM}.
By the first two chapters of \cite{FM} (now probably
the most difficult part of the story) this implies

\begin{Theorem}\label{ellmult}
Let $X \to C$ be an elliptic K\"ahler surface.	Then the multiplicities
of the elliptic fibration are determined by the underlying oriented
smooth manifold. In particular for K\"ahler elliptic surfaces
deformation equivalent, and oriented diffeomorphic are the same notions.
\end{Theorem}

This theorem has been well established with Donaldson theory by the work
of  Bauer, Kronheimer, Fintushel, Friedman, Morgan, Mrowka, O'Grady and
Stern.

\begin{Corollary}\label{plurigenera}
The plurigenera of a K\"ahler
surface are determined by the underlying oriented manifold.
\end{Corollary}

Let me remark that it seems to be known that
in the non K\"ahler case, with the exception of the equivalence of
deformation and diffeomorphism equivalence of non K\"ahler elliptic
surfaces, (where there can be a two to
one discrepancy) all the previous statements are true as well,
but seemingly for ``classical'' reasons like the homotopy type.

Inspired by results in	the preprint of Friedman and Morgan
I realised how the results in this article  give an easy proof of

\begin{Corollary}\label{poscurv}
No K\"ahler surface of non negative Kodaira dimension
admits a metric of positive scalar curvature
\end{Corollary}

While working on this article  a flood of information on the Seiberg
Witten classes came in. The holomorphic interpretation of the monopole
equations is already in Wittens paper \cite{Witten}, and it seems
that several people have remarked that his work implies that  the
canonical class is invariant for minimal surfaces of general type with
$p_g >0$ because of the numerical connectedness of the canonical
divisor. Kronheimer informed me that he, Fintushel, Mrowka,Stern and
Taubes are working on a note containing among many other things the
mentioned proof of the invariance of $\Kmin$.  The results and methods
of the before mentioned paper \cite{FM:SW} of  Friedman and Morgan are
rather similar to the present one. The main difference	seems to be that
they deal mostly with the case $p_g =0$,  and that they
rely on chamber changing formulas and a detailed analysis of the chamber
structure.   They also use a
 stronger version of the blow up formula  which allows them to	prove a
stronger version of theorem~\ref{main}.\ref{mainb}: if a surface of non
negative Kodaira dimension has a connected sum decomposition $X \iso  X'
\# N$,	where $N$ is negative definite, then $H_2(N,\Z) \subset
H_2(X,\Z)$ is spanned by $(-1)$-curves. We will indicate how this result
follow from the present methods.
Finally Taubes shows that the results
for K\"ahler surfaces are but the top of the iceberg. It seems that
most results can be generalised to symplectic manifolds
\cite{Taubes:sympl},\cite{Taubes:Gromov}.

\begin{Acknowledgment}
Thanks to  Stefan Bauer for pointing out a mistake in one of my
original arguments, and showing me the argument of
corollary~\ref{ellmult}. Thanks for Zhenbo Qin and Robert Friedman for
organising a very successful workshop in Stillwater, and for financial
support to attend. Thanks to Alexander Tichomirov and Andrej Tyurin  for
the opportunity to speak  at  the Yaroslav conference on algebraic
geometry on the  then possible invariance of $\Kmin$. Thanks to the
attendants of the Bielefeld Seiberg Witten Seminar (Stefan Bauer,
Manfred Lehn, Wei Ling, Viktor Pidstrigatch, Martin Schmoll, Stefan
Schr\"oer and Thomas Zink) for their comments and discussions.	Thanks
to Ian Hambleton for inviting me during April 1995 to the Max Planck
Institut f\"ur Mathematik in Bonn, the MPI is thanked for support.
Thanks to Steve Bradlow for pointing me to his and Garcia Pradas work on
the Vortex equation. Thanks to Hans Boden for not buying one of original
arguments on the localised Chern class. Finally special thanks to
Robert Friedman, for his insistent questions on the case $p_g =0$ which
revealed a serious mistake in my original exposition.
\end{Acknowledgment}

\section{Preparation}

We first prove the  corollaries from the main theorems~\ref{main}
and~\ref{ellmult}.

\begin{pf}
Corollary \ref{spheres}. Let $S$ be a positive sphere on a surface with
$\kod \ge 0$. Blow up $n =S^2 + 1$ times. Now $e = S + E_1 + \cdots
+E_n$ is a $(-1)$-sphere. Hence there is a $(-1)$-curve $E_0$ such that
$e = \pm E_0 \in H_2(X,\Z)$.  Then $S = \pm E_0$ , or $e = E_0 = E_1$
say.  The first possibility leads to the contradiction $E_0^2 \ge 0$,
the second to $S = 0 \in H_2(X,\Z)$. (Reducing non negative spheres to
$(-1)$-spheres is a well known trick, but I forgot where I read  it
precisely.)

Corollary~\ref{-infty} follows from corollary~\ref{spheres}.

Corollary~\ref{Kodaira}.By the above,
a K\"ahler surface is of Kodaira dimension $-\infty$ if it
contains a non trivial $(0)$-sphere.  Clearly  all ruled
surfaces contain one. To deal with $\P^2$, note that there is no surface
with $b_+ =b_1=0$. Thus diffeomorphisms between surfaces with $b_2=1$,
$b_1=0$ are automatically orientation preserving. Then a surface
diffeomorphic to $\P^2$ must contain a $(+1)$-sphere, and is therefore of
Kodaira dimension $-\infty$. Since $b_2 = 1$ it must in fact be equal to
$\P^2$ (alternatively use Yau's result that
$\P^2$ is the only surface with the homotopy type of $\P^2$
\cite[Theorem 1.1]{BPV}, but this is a
deep theorem).
We conclude that
Kodaira dimension $-\infty$ can be characterised by  just diffeomorphism
type. Without loss of generality we can therefore assume that
$\kappa \ge 0$.

If $\Kmin^2 > 0$, then $X$ is of general type. If $\Kmin^2 = 0 $
and $\Kmin$ is not torsion, then $\kod(X) = 1$, finally if $\Kmin$  is
torsion, $\kod(X) = 0$. This proves that Kodaira dimension is determined
by the oriented diffeomorphism type.  If $X$ and $Y$ are  orientation
reversing diffeomorphic, both are minimal, otherwise one of them would
contain a positive sphere. Then necessarily either  $K_X^2 = K_Y^2 = 0$,
or both have $K_X^2 ,K_Y^2 >0$, i.e. $X$ and $Y$ are of general type.
Now copy the argument of \cite[lemma S.4]{FM}: for minimal  surfaces
with $\kod = 0,1$, the signature  $\sign = \numfrac13(K^2 - 2e)\le 0$.
Thus  $\sigma(X) = -\sigma(Y) = 0$, and $e(X) = e(Y) = 0$. In Kodaira
dimension $0$, this leaves only tori and hyperelliptic surfaces, which
can fortunately be recognised by homotopy type \cite[lemma 2.7]{FM}.

Corollary~\ref{plurigenera}. Since $P_1=p_g$ is an oriented topological
invariant we will whence assume that $n\ge 2$.
We have to distinguish between the
different Kodaira dimensions.
For surfaces of general type (i.e $\kod = 2$) we argue as follows.
The plurigenera $P_n$ and $\chi(O_X)$ are birational invariants.
Then by Ramanujan vanishing and Riemann Roch (cf. \cite[corollary
VII.5.6]{BPV}) we have
\begin{equation}\label{gtpluri}
       P_n(X) = P_n(\Xmin) = \half n(n-1) \Kmin^2 + \chi(\O_X)
\end{equation}
 Since	$\chi(\O_X)$ is an  oriented topological
invariant the $P_n$ are oriented diffeomorphism
invariants in this case.
For surfaces with Kodaira dimension $0$ or
$1$ with a fundamental group that is not finite cyclic, we simply quote
\cite[S.7]{FM}. For surfaces with finite cyclic fundamental group, it
follows from the invariance of the multiplicities and the canonical
bundle formula which gives an explicit formula for $P_n(X)$ in terms of
the multiplicities and $\chi(O_X)$.  (see \cite[lemma I.3.18, prop.
I.3.22]{FM}).
\end{pf}

Here is an other easy corollary

\begin{Corollary}
Every $(-2)$-sphere $\tau$ is orthogonal to $\Kmin$ . If  there is a
$(-1)$-curve $E_1$ such that $\tau\cdot E_1 \ne 0$, then there is a
$(-1)$-curve $E_2$ such that $\tau = \pm E_1 \pm E_2 \in H_2(X,\Z)$.
\end{Corollary}

\begin{pf}
Let $R_\tau $ be the reflection in $\tau$. It is represented by a
diffeomorphism with support in a neighborhood of $\tau$.  By the
invariance of $\Kmin$ up to sign,  $R_\tau \Kmin = \Kmin + (\tau \cdot
\Kmin) \tau = \pm \Kmin$. But if $\Kmin \ne 0 \in H^2(X,\Q)$, then
$\tau$ and $\Kmin$ are indepent, since $\tau^2 = -2 $ and $\Kmin^2 \ge
0$. Thus in either case $(\tau, \Kmin) = 0$. Moreover if $E_1$ is a
$(-1)$-curve then  either $R_\tau E_1 = E_1$, $R_\tau E_1 = -E_1$,  or
there is a  different $(-1)$-curve $E_2$ such that  $R_\tau(E_1) = \pm
E_2$. The first possibility gives  $\tau\cdot E_1 = 0$, the second
$(\tau \cdot E_1)^2 = 2$ i.e. is impossible,  and the third  $(\tau\cdot
E_1) = \pm 1$. The statement follows.
\end{pf}

It will be convenient to first prove the main theorem~\ref{main}  with
(co)homology groups with $\Q$ coefficients,  and later mop up to prove
the theorem over $\Z$. Theorem~\ref{main} mod torsion is a formal
consequence of the existence of a set of basic classes
 $$
       \K(X) = \{K_1, K_2 \ldots \} \subset H^2(X,\Z)
 $$
functorial under oriented diffeomorphism between  $4$-manifolds with
$b_+ \ge 1$, and having the following properties:

\begin{Properties} \label{*}
For every K\"ahler surface $X$ of non negative Kodaira dimension
\begin{enumerate}
\item \label{i} the $K_i$ are of type $(1,1)$ i.e. represented by divisors,
\item \label{ii} if $X$ is minimal, then for every K\"ahler form $\Phi$,
	     $\deg_\Phi(K_X )\ge |\deg_\Phi(K_i)|$,
\item \label{iii} if $\~X \map{\sigma} X$ is the blow-up of a point $x
\in X$,
then $\sigma_*\K(\~X) = \K(X)$.
\item \label{iv} every $K_i$ is characteristic i.e.
	   $K_i \equiv w_2(X) \pmod 2$,
\item \label{v} $K_X \in \K$.
\end{enumerate}
\end{Properties}

In the case that $X$ is an algebraic surface  we could replace
item~\ref{ii} by weaker and
more geometric requirement that
$2g(H) - 2 \ge H^2 + |K_i\cdot H|$ for every
very ample divisor $H$ without changing the results.
We will see later that
%
%
Seiberg Witten theory will give such an inequality for all
surfaces minimal or not.
This should not be confused with a Thom conjecture type of
statement, since our methods do not give information about  the minimal
genus for  arbitrary smooth real surfaces in a homology class.
It is also clearly impossible to have a degree inequality like
property \ref{ii}
for all K\"ahler forms if $X$ is rational or ruled.

Recall that for algebraic surfaces, the Mori cone $\NEbar(X) \subset
H_2(X,\R)$ is the closure of the cone generated by effective curves. It
is dual to the nef (or K\"ahler) cone.	In other words, the numerical
equivalence class of a	curve $D$ lies in $\NEbar(X)$ if and only if
$H\cdot D \ge  0$ for all $H$ ample.  For a K\"ahler surface $(X,\Phi)$,
it will be convenient to define the
nef cone as closure of the positive cone in $H^{1,1}(X)
\subset H^2(X,\R)$ spanned by all K\"ahler forms, and  containing
$\Phi$.  The Mori cone $\NEbar$ is then just the dual cone in
$H_2(X,\R)\cap {H^{1,1}}^\dual$ i.e.
$$
       \NEbar=
 \{C \in  {H^{1,1}}^\dual \subset H_2(X,\R) \mid \int_C \omega \ge
    0, \txt{for all K\"ahler forms $\omega$}\}.
$$
(With this definition, a line bundle is nef iff for all $\epsilon > 0$,
it admits a metric such that the curvature form $F$ has
$\numfrac{\sqrt{-1}}{2\pi} F \ge -\epsilon \Phi$. A class $\omega \in
\NEbar$ if there exists a sequence of closed positive currents of type
$(1,1)$ converging to the dual of $\omega$, i.e is  $\NEbar$ dual to
$N_{\text{psef}}$ in \cite[proposition 6.6]{Demailly}. I am grateful to
Demailly for explaining this to me). We will freely identify homology
and cohomology by Poincar\'e duality.

\begin{Lemma} \label{decomp}
If a class $L \in H^{1,1}(X)$ satisfies $\deg_\Phi(K_X)\ge
|\deg_\Phi(L)|$ for all K\"ahler forms $\Phi$, then
there is a unique decomposition of the canonical divisor
$K_X = D_+ + D_- $ with $D_+$, $D_- \in \NEbar(X)$ such that
$L= D_+ - D_-$.
\end{Lemma}

\begin{pf}
Define $D_\pm = \half(K_X \pm L)$. Then  $K_X = D_+ + D_-$,
$L= D_+  - D_- $, and  $D_\pm \in \NEbar$.
\end{pf}

The following simple lemma is a minor  generalisation of the
fact that the canonical divisor of a surface of general type is
numerically connected \cite[VII.6.1]{BPV}.

\begin{Lemma} \label{connectedness}
Let $X$ be a minimal K\"ahler surface of non negative Kodaira
dimension. Suppose there is a decomposition $ K_X = D_+ + D_-$ with
$D_+$, $D_- \in \NEbar(X)\subset H^{11}(X)$. Then $D_+ \cdot D_- \ge 0
$, with equality if and only if say $K_X \cdot D_+ = D_+^2 = 0$. Thus
if $X$ is of general type then	$D_+ = 0$, if $\kod(X) = 1$, then  $D_+
= \lambda K_X$ with $0 \le \lambda \le 1$, and finally $D_+ = D_- = 0$
if $\kod(X) = 0$.
\end{Lemma}

\begin{pf}
First assume that $D_+^2 \le 0$.  Since $K_X$ is nef,  $D_+ \cdot D_- =
(K_X-D_+)\cdot D_+ \ge -D_+^2 \ge 0$, with equality iff $K_X \cdot D_+ =
D_+^2 = 0$. If $D_+^2 >0$  and $D_-^2 >0 $, then using the  K\"ahler
form $\Phi$, we can write $D_+ = \alpha \Phi + C_+$ and $D_- = \beta
\Phi + C_-$ with $\alpha$, $\beta >0$ and $C_\pm \in \Phi^\perp$.  By
the Hodge index theorem,
$$
       D_+ \cdot D_- = \alpha \beta \Phi^2 + C_+\cdot C_-
	     \ge \alpha\beta \Phi^2 - \sqrt{-C_+^2}\sqrt{-C_-^2} >0.
$$
The statement  for surfaces of general type follows directly
from Hodge index  and the fact that $K_X^2 > 0$.
If $\kod(X) = 1$, then	$K_X$ is
a generator of the unique isotropic subspace of $K_X^\perp$,
so $D_+ = \lambda K_X$, and $D_- = (1-\lambda) K_X$.
Since $K_X$, $D_+$ and	$D_- \in \NEbar(X)$, $\lambda$ is bounded by
$0 \le \lambda \le 1$.
Finally if $\kod(X) = 0$, $K_X$ is numerically trivial and,
$D_+$ and $D_-$ must be zero as well.
\end{pf}

\begin{Lemma} \label{inequality}
Let  $X$ be a surface of non negative Kodaira dimension with
$(-1)$-curves $E_1, \ldots E_m$. Assume that $\K$ has
properties~\ref{*}.  Then $K_i^2 \le K_X^2$ for all $K_i \in \K$, with
equality if and only if
$$
       K_i = \lambda \Kmin + \sum \pm E_i \in H^2(X,\Q)
$$
where $\lambda = \pm 1$ if $X$ is of general type,  $\lambda$ is a
rational number with  $|\lambda | \le 1 $ if $\kod(X) = 1$,	and
where $\lambda = 0$ if $\kod(X) = 0$.
\end{Lemma}

\begin{pf}
By property  \eqref{iii}, and \eqref{iv},
$K_i = K_{i,\min} + \sum_j (2a_{ij} + 1)E_j$.
Thus
$$
       K_i^2  \le K_{i,\min}^2 - \#(-1)\hbox{-curves},
$$
with equality if and only if $a_{ij} =0$, or  $-1$ for all $i,j$.  Since
$K_X^2 = \Kmin^2 - \#(-1)$-curves, we can assume that $X$ is minimal.
Using property~\eqref{i} and \eqref{ii} and lemma~\ref{decomp}, write
$K_X = D_+ + D_-$ and $K_i = D_+ - D_-$, with $D_\pm \in \NEbar(X)$.
Then by lemma~\ref{connectedness} $K_i^2 = K_X^2 - 4 D_+\cdot D_- \le
K_X^2 $ with equality under the stated condition. Note that this lemma
does not use diffeomorphism invariance, nor that $K_X \in \K$.
\end{pf}

We are now in a position to formulate and  prove  half of the main theorem

\begin{Proposition} \label{Kcharacterisation}
Assume that for all 4-manifolds $X$ with $b_+ \ge 1$ there is a set of
basic classes  $\K(X) = \{K_1, K_2, \ldots\} \subset H^2(X,\Z)$
functorial under oriented diffeomorphism having properties~\ref{*}.
Then $\Kmin$ is an oriented $C^\infty$ invariant up to sign and
torsion, and every $(-1)$-sphere is represented by  a $(-1)$-curve up
to sign and torsion.
\end{Proposition}

\begin{pf}

Using  lemma~\ref{inequality} we can easily reduce the invariance of
$\Kmin$ up to sign and torsion to
showing  that $(-1)$-spheres are represented by $(-1)$-curves up to sign
and torsion.

Since $K_X \in \K$, there is  a nonempty subset
$\K_0  = \{ K_j\} \subset \K$ with
$K_j^2 = K_X^2 = 2e(X) + 3\sign(X)$.
Consider the projection $K_{j,\min}$ of $K_j$ to the minimal model i.e.
the projection to the orthogonal complement of the $(-1)$-spheres.

If $K_{j,\min}^2 >0$, then by lemma~\ref{inequality}, $X$ is of general
type, and $K_{j,\min} = \pm \Kmin$ up to torsion.
If $K_{j,\min}^2 = 0$, there are two possibilities. If $K_{j,\min}$ is
torsion for all $j$, then again by lemma~\ref{inequality},
$X$ is of Kodaira dimension $0$ i.e. $\Kmin$ is also
torsion.
Otherwise we choose $j$ such that $K_{j,\min} \ne 0$ has maximal
divisibility. Since $K_X \in \K_0$ our little lemma shows that,
the Kodaira dimension is $1$ and $K_{j,\min} = \pm \Kmin$.

Now let $e$ be the class of a $(-1)$-sphere in $H^2(X,\Q)$. Without loss of
generality, we can assume that $K_X\cdot e < 0$.
Consider $R_e$ the reflection generated by a $(-1)$-sphere $e$. It is
represented by an orientation preserving  diffeomorphism.
Since $\K$ is invariant under oriented diffeomorphisms, the
characterisation of basic classes with square $K_X^2$ tells us that
\begin{align}
	R_e K_X &=  \Kmin  +  \sum  E_i  + 2 (K_X\cdot e)e
\label{line1}
\\
		    &=	\lambda \Kmin + \sum \pm E_i
\label{line2}
\end{align}
with $|\lambda| \le 1$.
Taking	intersection with $E_i$ we find that $(E_i\cdot e)( e\cdot K_X)
=0$~or~$1$. Since $ K_X\cdot e \equiv e^2$ is odd, $e$ is either
orthogonal to all $(-1)$ curves (i.e. $ e\in H^2(\Xmin,\Q)$) or there is
a $(-1)$-curve, say $E_1$, such that
$K_X\cdot e =E_1 \cdot e= -1$. However, $e \in H^2(\Xmin)$ implies
that $e = {\lambda -1 \over 2K_X \cdot e}\Kmin$,
which is impossible because $\Kmin^2 \ge 0$.
Thus, after renumbering the $(-1)$-curves, \eqref{line1} and
\eqref{line2} can be rewritten to
\begin{equation}\label{(-1)-sphere}
       e = \half(1-\lambda) \Kmin + \sum _{i=1}^N E_i
\end{equation}
with $N = \quart(1-\lambda)^2 \Kmin^2 + 1$.

Now reflect $e$ in $E_1^\perp$. $R_{E_1} e$ is yet another
$(-1)$-sphere,
 so it has a  representation as in equation~\eqref{(-1)-sphere}, except
possibly for an overall sign
\begin{align*}
       R_{E_1}e &= \half(1-\lambda)\Kmin - E_1 + \sum _{i=2}^N E_i
\\
	      &= \pm \big(\half(1-\mu) \Kmin + \sum_{j=1}^M E_{i_j}\big).
\end{align*}
Upon comparison, we see that the sign is minus, that $N=M =1$, and that
$0 \le 1-\lambda = \mu - 1 \le 0$ unless
$\Kmin = 0$. In  other words $ e = E_1 \in H^2(X,\Q)$.
\end{pf}

\section{The  localised Euler class of a Banach bundle.}\label{top}

We will use a construction pioneered by Pidstrigatch and Pidstrigatch
Tjurin \cite{Pidstrigatch:instanton}, \cite[\S 2]{PT}, which is a
convenient and general way to define fundamental cycles for moduli
spaces arising from elliptic equations.  Unfortunately their
construction is not quite in the generality we will need it, and we will
therefore set it up in fairly large generality here. The cycle is the
localised homological Euler class of an infinite dimensional bundle. It
can be used to give definitions that avoid transversality arguments
needing small deformations, generic metrics etcetera, although
transversality will be extremely useful for computations and proofs. The
construction is modeled on Fultons intersection theory and in the
complex case it makes the whole machinery of excess intersection theory
available. However, although the construction is very simple in principle,
the whole thing has turned  a bit technical. On
first reading it is best to ignore the difference between \v Cech and
singular homology, and continue to proposition \ref{locEuler}, the
construction of the Euler class in the proof of this proposition
and corollary~\ref{locChern}.  Some readers might even
want to continue to the next section, since we will use rather little of
the general machinery for the proofs of the theorems and corollaries in
the  introduction.

\smallskip\noindent
We first make some algebraic topological preparations.
For any pair of topological spaces $A \subset X$,
homology with closed support and
with local coefficients $\locsys$ is defined as
$$
       H_i^{\cl}(X,A;\locsys) = \invlim_K H_i(X,A \union (X-K);\locsys)
$$
where we take the limit over all compacta $K \subset X - \o A$.
$H^{\cl}_*$ is functorial under proper maps.  Unfortunately this
``homology theory'' suffers the same tautness problems that singular
homology has. To be able to work with well behaved cap products we will
have to complete it.  The following works well enough for our purposes but
is a bit clumsy.

Suppose that $X$ is {\sl locally modelable} i.e.
is locally compact Hausdorff and has local models
which are each subsets of some $\R^n$. Obviously locally compact
subsets of locally modelable spaces are locally
modelable, in particular a closed subset of a local modelable space is
locally modelable.
Then for every compact subset $K \subset X - \o  A$ there is a
neighborhood $U_K \supset K$ in $X$ which embeds in $\R^N$.
We now define
$$
       \cH_i^{\cl}(X,A,\locsys) = \invlim_K \cH_i(U_K,A\cap U_K \union
(U_K -K); \locsys)
$$
where for every pair $(Y,B)$ in a manifold $M$, \v Cech homology is defined
as
$$
       \cH_i(Y,B) = \invlim \{H_i(V,W),\  (V,W) \txt{neighborhoods of}
(Y,B) \txt{in} M\}
$$
This definition depends neither on  the choice of $U_K$,
nor on the embedding $U_K \inj\to \R^N$, since two embeddings are
dominated by the diagonal embedding, and $\cH_*(Y,B)$ does not depend on
$M$ but only on $(Y,B)$ (c.f \cite[VIII.13.16]{Dold}).

Fortunately we do not usually have to bother with \v Cech homology.
Suppose in addition that $X$ is locally contractible e.g. locally a
sub analytic set (c.f. \cite[\S I.1.7]{GoreskyMcPherson}, and the fact
that Whitney stratified spaces admit a triangulation).	Then $X$ is
locally an Euclidean neighborhood retract (ENR) by   \cite[IV
8.12]{Dold} and since in a Hausdorff space a finite union of ENR's  is
an ENR by \cite[IV 8.10]{Dold} we can assume that $U_K$ is an ENR. Now
assume that $A$ is open. Then by \cite[prop. VIII 13.17]{Dold}
$$
     \cH_* (U_K,U_K\cap A \union (U_K-K)) \iso
      H_*(U_K,U_K \cap A \union (U_K - K)) \iso H_i(X,A \union X-K).
$$
Thus in this case $\cH_*^{\cl}(X,A) = H_*^{\cl}(X,A)$. If $A$ is closed
and locally contractible then one should be able to organise things such
that $U_K \cap A$ is an ENR and the same conclusion would hold.

\begin{Lemma}\label{capproduct}
Let $X$ be a locally modelable space, and $Z$ a locally compact (e.g.
closed) subspace, then there are cap products
$$
       \cH^i(X,X - Z,\locsys) \tensor \cH^{\cl}_j(X,\locsys') \map{\cap}
	       \cH^{\cl}_{j-i}(Z, \locsys\tensor\locsys')
$$
with the following properties.
\begin{enumerate}
\item
If $Y$ is locally embeddable, $f: Y\to X$ is proper, and
$\sigma' \in \cH^{\cl}_j(Y, Y-f^{-1}(Z))$, then the push-pull formula
holds:
$$
       f_*(f^* c \cap \sigma') = c \cap f_* \sigma'.
$$
\item
If $Z \inj\map{i} W$ is proper and $W$ is locally compact,
we can increase supports i.e.
$$
	c|_{(X,X-W)} \cap \sigma = i_*(c \cap \sigma).
$$
\end{enumerate}
\end{Lemma}

\begin{pf}
For every $c \in \cH^i(X,X-Z)$ and $\sigma \in \cH_j^{\cl}(X)$,
we have to construct a class $c \cap \sigma \in \cH_{i-j}(Z, Z-K)$
for a cofinal family of compacta $\{K\}$. Since  $Z$ is locally compact,
every compactum $K$ is contained in a  compactum $L
\subset Z$ with $L \supset \o L \supset K$. Likewise there exists a
compactum $L' \supsupset L$.
By excision it suffices to construct a class in $\cH_{i-j}(L,L-K)$. Let
$U_{L'}$ be a neighborhood of $L'$ in $X$ which embeds in $\R^N$.
Let $V_L$, $W_{L-K} \subset V_L$, and $V_K \subset V_K$ be neighborhoods
of respectively $L$, $L_K$ and $K$  in $\R^N$. Define $U_L = V_L \cap
U_{L'}$.  We can assume that $U_L \cap Z = U_L \cap L'$, $V_K \cap Z =
V_K \cap L$, and after replacing $V_{L-K}$ by $(V_{L-K} - (L' \cap
W_{L-K}^c) \union V_K$, that $V_L \cap (L'-K) = W_{L-K} \cap (L' - K)$.
Then our task is to construct a class $c_L\cap \sigma_L \in
H_{i-j}(V_L,W_{L-K})$ possibly after shrinking $V_L$ and $W_{L-K}$.

We have a restriction map $\cH^i(X,X-Z) \to \cH^i(U_L,U_L-L')$. After
shrinking $V_L$ if necessary, $c|_{(U_L,U_L-L')}$ comes from a
class $c_L \in H^i(V_L,V_L -L')$.
By definition there is map
$$
       \cH^{\cl}_j(X) \to \cH_j(U_L,U_L-K) \to H_j(V_L,V_L-K).
$$
Let $\sigma_L \in H_j(V_L,V_L-K)$ be the image of $\sigma$.

Now write $V_L- K = (V_L - L') \union (W_{L-K} - K)$. Then the standard
cap product \cite[VII Def. 12.1]{Dold} gives a map
$$
	H^i(V_L,V_L-L') \tensor H_j(V_L,V_L-K) \map{\cap}
H_{j-i}(V_L,W_{L-K}-K)
$$
so we get a class $c_L \cap \sigma_L \in H_{j-i}(V_L, W_{L-K})$ as
required. Since if $K' \supset K$, choices for $K'$ will work a
fortiori for $K$, we can pass to the limit.

To prove the first property, note that since $f$ is proper, $f^{-1} Z$
is locally compact. Choose compacta $K \subsubset L \subsubset L'
\subset Z$ giving compacta $f^{-1}K \subsubset	f^{-1}L \subsubset
f^{-1}L'$.  Note that compacta of the form $f^{-1}K$ are a cofinal
family of compacta in  $f^{-1}(Z)$. Embed neighborhoods $U_{L'} \subset
V_{L'} \subset \R^N$  and $U_{f^{-1}L'} \subset \R^M$.	Now we
carry out the construction above  with the diagonal embedding of
$U_{f^{-1}L'}$ in  $ \R^{N + M}$. Let $V_{f^{-1}L'}$ be a neighborhood
of $U_{f^{-1}L'} \in \R^{N+M}$.  We can assume that $V_{f^{-1}L'} \to
V_{L'}$ under the projection $\pi$ to $\R^N$. We can also  assume that
$c|_{(U_L,U_L-L')}$ comes from a class $c_L \in H^i(V_L,V_L - L')$.
Finally let $\sigma_{f^{-1}L'}$ be an image of $\sigma$ in
$H_j(V_{f^{-1}L'}, \pi^{-1}W_{K-L})$.  Then the first property follows
from the identity
$$
       \pi_*(\pi^* c_L \cap \sigma'_{f^{-1}L'}) = c_L \cap \pi_*
       \sigma'_{f^{-1}L'}
$$
in $H_j(V_l, W_{K-L})$. The second property is left to	reader.
\end{pf}

A smooth manifold $X$ of dimension $n$, has an orientation system
$\orr(X)$, the sheafification of the presheaf  $U \to H^n(X,X-U)$.
Equivalently, we can define  $\orr(X)$ as the sheaf
$R^d\pi_*(X\times X, X\times X
-\Delta, \Z)$, where $\Delta$ is the diagonal of $X\times X$, $\pi$
the projection on the first coordinate, and $R^d\pi_*$ the
parametrised version of the $d$ th cohomology.

Likewise
for a real vector bundle $E$ of rank $r$ there is an orientation system
$\orr(E)$, the sheafification of $H_q(E|_U,E|_U - U)$.
We have $\orr(X) = \orr (TX)^\dual$,
as can be seen immediately from the alternative description of $\orr(X)$
and excision.

A manifold $X$ has a unique fundamental class
$[X] \in H^{cl}_n(X,\orr(X))$  in singular or \v Cech homology such that
for small $U$
$$
    [X]|_{X-U}\in H_d(X,X-U,H^d(X,X-U))= \Hom(H^d(X,X-U),H^d(X,X-U))
$$
is identified with the identity (cf \cite[p. 357]{Spanier}).

Similarly, a bundle $E$  has a
Thom class $\Phi_E \in \cH^r(E,E-X, \orr(E))$ \cite[p. 283]{Spanier}.
In turn for every section $s$ in $E$ with zero set $Z(s)$,
the Thom class defines a localised cohomological
Euler class $e(E,s) = s^*\Phi_E \in \cH^r(X, X - Z(s), \orr(E))$.

\smallskip\noindent
Let $M$ be a Banach manifold, $E$ a real Banach vector bundle on $M$
and $s$ a section of $E$ with zero set $Z(s)$.
The section induces an exact  sequence
\begin{equation}\label{tosplit}
       0 \map{} E \map{} s^* TE \map{\pi} TM \to 0,
\end{equation}
which expresses that the vertical tangent bundle of the total space of
$E$ is canonically  isomorphic to the bundle $E$. On $Z(s)$ we have a
canonical splitting of this sequence,
given by the sequence
$$
       0 \map{} TM \map{Ts_0} s_0^*TE \map{} E \map{} 0
$$
defined by the zero section $s_0$, and the
identification $s^*TE|_{Z(s)} = s_0^*TE|_{Z(s)}$ over $Z(s)$.
This gives a canonical map
$$
       Ds: TM|_{Z(s)} \map{Ts} s^*TE = s_0^*TE \map{} E|_{Z(s)}.
$$

If $D$ is a connection on $E$ then
$D(s)$ is a splitting that extends the canonical splitting over $Z(S)$
(hence the notation) but in general connections need not exist on Banach
manifolds. We will avoid choosing non canonical splittings.


\begin{Proposition}\label{locEuler}
Let $M$ be a smooth Banach manifold, $E$ a banach bundle over $M$ and $s$
a section in $E$. Assume that
\begin{enumerate}
\item
{\sloppy
The map $Ds$ is a section in the bundle $\Fred^d(TM|_{Z(s)},E|_{Z(s)})$
of Fredholm maps of index $d$. We say that
$Z(s)$ has virtual dimension $d$, and that $Ds$ is Fredholm of index $d$.
}
\item
The real line bundle
$\det(\Ind(Ds))$ is trivialised over $Z(s)$.
\end{enumerate}
Then these data define a \v Cech homology class with closed support
$$
       \Z(E,s)= \Z(s) \in \cH^{\cl}_d(Z(s),\Z)
$$
with the following properties.
\begin{enumerate}
\item\label{smoothcase}
The class $\Z(s) = [Z(s)]$ if $Z(s)$ is smooth of dimension $d$ and
carries the natural orientation defined by the trivialisation of
$\det(\Ind Ds)$,
\item\label{homotopy}
if $\{C\}$ is a family of closed subsets of $M$ such that $C \cap Z(s)$ is
compact for all $C$, then there is a natural map $\cH_j(Z(s)) \to
\invlim_C H_j(M,M-C,\Z)$, and if $s_t$ is a one parameter family of
sections with this property then
$\Z(s_0) = \Z(s_1)) \in \invlim_C H_d(M,M-C,\Z)$.
\end{enumerate}
For every exact sequence
$$
       0 \to E' \to E \to E'' \to 0,
$$
defined over a neighborhood of $Z(s)$,
let $s''$ be the induced section in $E''$, and $s'$ the induced
section of $E'|_{Z(s'')}$ with zero set $Z(s)$. Then
\begin{enumerate}
\setcounter{enumi}{2}
\item \label{fdeuler}
if $E'$ has finite rank
$$
       \Z(s) = e(E'|_{Z(s'')},s')\cap \Z(s''),
$$
\item \label{stability}
if $Ds''|_{Z(s)}$ is surjective,
then $Z(s'')$ is smooth in a neighborhood of $Z(s)$,
$Ds':TZ(s'')|_{Z(s)} \to E'|_{Z(s)}$
is Fredholm, with $\Ind Ds' \iso \Ind Ds$   and
$$
       \Z(E,s) = \Z(E'|_{Z(s'')},s').
$$
\end{enumerate}
\end{Proposition}

For property \ref{homotopy} there are two typical situations we have in
mind. One is that we have a natural connected family of sections $s_t$
such that $Z(s_t)$ is compact for all $t\in T$. In this situation we get
a   homology class $\Z(s_{t_0}) \in H_d(M)$ independent of the choice of
$t_0$ (take $\{C\} = \{M\}$). Such will be  the case in Seiberg Witten
theory.  In the other case we again have a family of sections $s_t$ but
there is ``bubbling'' which invariably means we lack some a priori
estimate. For example in Donaldson theory, the moduli space of ASD
connections with curvature bounded in the $L^4$ norm is compact.
Therefore  it is natural to define a   family of subsets $\{\B^{\le
C}\}_{C \in \R^+}$  in the space $\B^*$ of all irreducible $L^2_2$
connections mod gauge, where $\B^{\le C}$ the subset of  connections
with $L^4$ norm of the curvature bounded by $C$.

\begin{pf}
If $M$ (hence $E$) is a finite dimensional manifold of dimension $N+d$
then  $E$ is a real  vector bundle of rank $N$ with an isomorphism
$\det(E) = \det(TM)$ over $Z(s)$.
Let $[M] \in H^{\cl}_N(M,\orr(M))$ be the fundamental class,
and $\Phi_E$ the twisted Thom class of $E$ in
$H^{N-d}(E,E-M,\orr(E))$.
Define
$$
       \Z(s) = e(E,s) \cap [M] \in \cH^{\cl}_d(Z, \orr(E)\tensor \orr(M))
	   = \cH^{\cl}_d(Z(s),\Z)
$$
i.e. $\Z(s)$ is the Poincar\'e dual of the localised  cohomological
Euler class.   In the last step we used the chosen trivialisation of
$\orr(E) \tensor \orr(M) = \orr(\det TM^\dual\tensor\det E) =\orr
(\det(\Ind(Ds)))$  given by the trivialisation of the index.

In the infinite dimensional case we proceed similarly but we have to go
through a limiting process and use that we know what to do when the
section is regular.  For each
compactum $K \subset Z$ we have to construct a class $\Z_K \in
\cH_d(Z,Z-K)$ such that for  $K' \supset K$ the class $\Z_{K'}|_{Z-K} =
\Z_K$  under the restriction map $H_d(Z,Z-K') \to
H_d(Z,Z-K)$.

Over  a neighborhood $U$ of $K$ in $M$ we can find a
finite rank $N$ subbundle $F$ of $E$ such that $\Im(Ds)|_K  + F|_K =
E|_K$. Such a bundle certainly exists: we can choose a finite number of
sections $s_1, \ldots s_N$ such that the $s_i$ span
$\Coker(Ds_x)$ for
every $x \in K$, and possibly after perturbing we can assume that the
$s_i$ are linearly independent in a neighborhood.
Let $\~E$ be the quotient bundle $E/F$ defined over $U$, and $\~s$ the
induced section with zero set $M_f = Z(\~s)$  ($f$ is for finite, $M$ is
for, well, manifold).

Clearly the map $TM|_{Z(s)} \map{Ds} E|_{Z(s)} \map{} \~E $ is
surjective.  Since the canonical map $D\~s$ on $M_f$ restricts to this
composition on $Z(s)$,	$D\~s$ is surjective on $M_f$ possibly after
shrinking $U$.	Hence $M_f$ is a smooth manifold.  Let $T =
\ker(TM|_{M_f} \to \~E)$.  There is a canonical identification $T \iso
TM_f$. Now $T$ is a bundle of rank $N+d$ since
\begin{equation}\label{theindex}
       \Ind(Ds)|_K = T - F.
\end{equation}
Thus  $M_f$ has dimension $N+d$.

On $M_f$, the section $s$ in $E$ lifts to a section $s_f$ of the
subbundle $F$. Clearly $Z(s_f) = Z(s) \cap U$.
Define
$$
   \Z_K  = e(F|_{M_f},s_f) \cap [M_f] \in \cH_d(Z(s), Z(s) - K;\Z)).
$$
Here we have used the restriction map
$$
       \cH^{\cl}_d(Z(s)\cap U; \orr(F)\tensor \orr(M_f))
		   \to\cH_d(Z(s),Z(s)-K; \orr(F) \tensor \orr(M_f)),
$$
the identification
$\orr(\det(\Ind(Ds))) = \orr(F) \tensor \orr(M_f)$
and the chosen trivialisation of $\det(\Ind(Ds))$
as in the finite dimensional case.

This construction does not depend on the choices. If $F$ and $F'$ are
two choices of subbundles of $E$ then there is third bundle
$F''$ containing $F + F'$. We can therefore assume that $F$ is a
subbundle of $F'$. Then using primes to denote objects we get out of the
construction above using $F'$ instead of $F$, we have a section
$s'_f$ in $F'$, a  section $s''_f$ in $F'/F$ cutting out $M_f$ in $M'_f$
and the identity
\begin{align*}
       \Z_K' &= e(F'|_{M'_f} s'_f)\cap [M'_f]
\\
	     &= e(F|_{M_f},s_f)\cap e(F'/F|_{M'_f}, s''_f) \cap [M'_f]
\\
	     &= e(F|_{M_f},s_f) \cap [M_f] = \Z_K.
\end{align*}
Note that in the third step we have used the identification
$\orr(M_f) = \orr(M'_f) \tensor \orr(F'/F))|_{M_f}$.
In particular, if $K' \supset K$ all choices on $K'$ work a fortiori for
$K$,  so we can pass to the limit.

The relation $\Z(s) = [Z(s)]$ for regular sections
(property~\ref{smoothcase}), and the
compatibility with Euler classes of finite rank bundles
(property~\ref{fdeuler})
are now clear from the construction.
The stability property~\ref{stability} also follows from
the construction. For every compactum $K$,
we can choose the finite rank subbundle $F$
as a subbundle of $E'$.  Then $\~E \surj\to E''$. Now one checks
that by a diagram chase that
$$
	  Z(\~E,\~s) = Z(E'/F|_{Z(E'',s'')}, s'\bmod F)
$$
and that
{
\let\to=\rightarrow
\begin{align*}
       TZ(\~E,\~s) &= \Ker(TM \to \~E)
\\
		    &= \Ker(\Ker(TM \to E'') \to E'/F)
\\
		    &= \Ker(TZ(s'') \to E'/F) = TZ(E'/F).
\end{align*}
}
In particular, the orientations agree. Thus we see that
\begin{align*}
       \Z_K(E,s) &= e(F,s_f)\cap [Z(\~E,\~s)]
\\
		  &= e(F,s_f)\cap [Z(E'/F|_{Z(E'',s'')}, s'\bmod F)]
		  = \Z_K(E'|_{Z(s'')},s').
\end{align*}
It only remains to pass to the limit over $K$.

To see that $\cH^{\cl}_j(Z(s))$ maps to $\invlim_C H_j(M,M-C)$ note that
for every compact subset $K = C \cap Z(s)$, we constructed a finite
dimensional manifold $M_f \supset K$.
Then we have maps
\begin{align*}
       \cH_j^{\cl}(Z(s)) \to &\cH_j(Z(s),Z(s)-K)
		 = \cH_j(Z(s) \cap M_f,  Z(s) \cap M_f -K)
\\
		&\to H_j(M_f, M_f - C) \to H_j(M,M-C).
\end{align*}
Again this map is independent of choices, and we can pass to the limit.

The homotopy property of $\Z$ is a formal consequence of the
compatibility with finite dimensional Euler classes. Consider the
trivial bundle $\R$ over the interval $[-1,2]$ with the one parameter
family of sections  $\theta - \tau$ where $\theta :[-1,2] \to \R$ is the
inclusion and $0 \le \tau \le 1$. Then clearly $e(\R,\theta) =
e(\R,\theta-1) \in H^1([-1,2],\{-1,2\})$ is the canonical generator.
Consider  $M\times [-1,2]$. Let $\pi:M	\times [-1,2] \to M$
be the projection and $S: M\times [-1,2] \to \pi^*E$
an extension of our one parameter family of sections e.g. $S_t = s_0$
for $t \le 0$ and  $S_t = s_1$ for $ t\ge 1$.
The bundle $\pi^*E \directsum \R$ has a one parameter family of
sections $(S,\theta- \tau)$. Now
\begin{align*}
       \Z(s_0) &\buildrel \ref{stability}\over=
		  \pi_* \Z(\pi^*E \directsum \R; (S,\theta))
\\
		&\buildrel \ref{fdeuler}\over=
		      \pi_* e(\R,\theta) \cap \Z(\pi^*E;S)
\\
		&= \pi_* e(\R,\theta - 1)\cap \Z(\pi^*E;S)
\\
		&= \pi_* \Z(\pi^*E \directsum \R;(S,\theta-1))
		= \Z(s_1)
\end{align*}
\end{pf}

\begin{Corollary}\label{locChern} (compare \cite[prop. III.2.4]{PT})
Let $M$ be a complex Banach manifold, $E$ a holomorphic vector bundle
and $s$ a holomorphic section with zero set $Z(s)$
Assume that $Ds$ is a section of $\Fred_\C^d(TM|_{Z(s)},E|_{Z(s)})$.
We say that $Z(s)$ has complex virtual dimension $d$, and that $Ds$ is
Fredholm of complex index $d$.
Then the localised Euler class $\Z(s)= [Z(s)] \in H^{\cl}_{2d}(Z(s),\Z)$,
if $Z(s)$ is a local complete intersection
of dimension $d$, and more generally
\begin{equation}\label{magic}
       \Z(s) = [c(\Ind(Ds))^{-1} c_*(Z(s))]_{2d}
\end{equation}
where $c_*(Z(s))$ is the total homological chern class of $Z(s)$ defined
analogous to \cite[example 4.2.6]{Fulton} by equation \eqref{homchernclass}
and coincides with the Poincar\'e dual of the
cohomological chern classes of the tangent bundle if $Z(s)$ is
smooth.
\end{Corollary}

\begin{Remark}
If $Z(s)$ is smooth we can
even get away with an almost complex manifold $M$ and the assumption
that  $Ds$ is complex linear.
\end{Remark}

\begin{Remark}
I have tacitly removed $M$ and $E$ from the notation of the
homological Chern class $c_*(Z(s))$. I strongly believe that
$c_*(Z(s))$ is independent of the embedding but I did not prove this.
There is one case where independence of $c_*(Z(s))$ on the embedding can
be proved completely analogous to \cite[Example 4.2.6]{Fulton} by
simply replacing algebraic arguments by complex analytic ones: if for
every $K \subset Z(s)$ compact, there exists a {\em holomorphic}
finite rank sub bundle $F \inj\to E$ defined over a neighborhood of
$K$ such that $F|_K + \Im(Ds)|_K = E|_K$. Then a neighborhood $U_K$ of
$K$ in $Z(s)$ sits in a complex rather then almost complex finite
dimensional manifold
$M_f$.	Such a bundle should typically exist if $Z(s)$ has the
structure of a quasi projective variety, and $\Coker DS$ has the
interpretation of a coherent sheaf as in \cite[\S 5,
\S6]{Pidstrigatch:instanton}.
%
\end{Remark}

\begin{pf}
We will use  Mac Phersons graph construction, that is we consider  the
limit  $\lambda \to \infty$ of the map $(\lambda s: 1)$ in $\P(E \oplus
\O)$ or finite dimensional approximations thereof.  We use the notations
of the proof of proposition~\ref{locEuler}.

For a compactum $K \subset Z(s)$ we choose the finite rank bundle $F$ as
follows. It is a  complex bundle, and in every point of $Z(s)$ there are
sections of $F$  which restricted to a neighborhood are  holomorphic
sections of  $E$ and which span locally  a subbundle $F^{\hol} \inj\to
F$, such that  $Ds: TE|_{Z(s)} \surj\to E/F^{\hol}|_{Z(s)}$ is a
surjection.    We do not assume that $F$ is a holomorphic subbundle,
because I do not see a reason why such a bundle should exist. However
since  $F$ is a complex bundle, the quotient  $\~E = E/F$ and
$$
       TM_f|_{Z(s)} = T|_{Z(s)} =
	\Ker (TM|_{Z(s)} \map{Ds} E|_{Z(s)} \to \~E|_{Z(s)})
$$
are complex bundles.
We extend this complex structure on $TM_f$
over all of $M_f$, possibly after shrinking $M_f$,
making it into an almost complex
manifold of complex dimension $d+ N$.

\def\0{\overline{(0,1)}}
\def\lsf{\overline{(\lambda s_f,1)}}
\def\sfnul{\overline{(s_f,0)}}
\def\sflinv{\overline{(s_f,1/\lambda)}}

Consider the space $\P(F \oplus \O) \map{\pi} M_f$. Then the total space
of $F$ embeds in  $\P(F \oplus \O)$. The image of the zero section will
also be called the zero section, and the complement of $F$ the divisor
at infinity.

Let $Q$ be the universal quotient bundle.
The bundle $Q$ has sections $\0$ , and $\lsf$
cutting out the zero section and the graph of $\lambda s_f$
respectively. Equivalently we can cut out the graph of $\lambda s_f$ by
$\sflinv$. Then clearly as $\lambda \to \infty$ the graph degenerates to
a set contained in the zero set of $\sfnul$.

Now $Z(\sfnul)$ has two ``irreducible components''. One
component $\~M_f \inj\to \P F$ is the closure of the image of
$(s_f:0): M_f -Z(s) \to \P F|_{M_f - Z(s)} \subset \P(F\oplus \O)$
It will be called the strict transform.
The other component is just  $\P(F \oplus \O)|_{Z(s)}$.   Let $\E_f =
\~M_f \cap \P(F \oplus \O)|_{Z(s)}$. It will be called the exceptional
divisor.

I claim that
\begin{equation}\label{dimclaim}
       \cH^{\cl}_{2d+2N-1 + i}(\E_f) =0 \txt{for} i \ge 0.
\end{equation}
Accepting this claim we see that
$\~M_f$ carries a unique fundamental class $[\~M_f]$
restricting to $[\~M_f - \E_f]$ by the exact sequence
$$
       \cH_{2d + 2N}^{\cl}(\E_f) \to
	\cH_{2d+2N}^{\cl}(\~M_f) \to
	H^{\cl}_{2d+2N}(\~M_f - \E_f) \to
	\cH^{\cl}_{2d+2N - 1}(\E_f)
$$
Consider $ C' = \Z(\sfnul) - [\~M_f] \in \cH^{\cl}_{2d +
2N}(Z(\sfnul))$.  Then $C'$ comes from a unique class $C \in
\cH^{\cl}_{2d+2N}(\P(F \oplus \O)|_{Z(s)})$ because of the sequence.
$$
       0 \to \cH^{\cl}_{2d+2N}(\P(F \oplus \O)|_{Z(s)}) \to
	 \cH^{\cl}_{2d+2N}(Z(\sfnul)) \to H^\cl_{2d+2N}(\~M_f - \E_f).
$$


Now note that $Q$ restricted
to the zero section is canonically isomorphic to $F$.
We therefore have the following chain of equivalences
\begin{align*}
       \Z(s)_K &= e(F,s_f) \cap [M_f]
\\
	     &=\pi_* e(\pi^*F,\lambda s_f) \cap e(Q,\0) \cap [\P(F\oplus\O)]
\\
	     &=\pi_* \(e(Q,\lsf)\cup e(Q,\0)\) \cap [\P(F \oplus \O)]
\\
	     &=\pi_* e(Q,\0) \cap e(Q,\sflinv) \cap [\P(F \oplus \O)]
\\
	     &=\pi_* e(Q,\0) \cap \(e(Q,\sfnul)\cap [\P(F \oplus\O)]\).
\\
	     &=\pi_* e(Q,\0) \cap \Z(\sfnul)
\end{align*}
If we accept the claim~\eqref{dimclaim}
for a moment and we note that the support of
$\~M_f$ and $e(Q,\0)$
are disjoint we see further that
$$
       \Z(s)_K = \pi_*e(Q,\0) \cap C' = \pi_* e(Q) \cap C
$$
If we use that	$e(Q)= c_\Top(Q)$
 this can be rewritten further to
\begin{align*}
	\Z(s)_K &= [\pi_* c(Q)\cap C)]_{2d}
\\
	      &= [c(F) \pi_*\((1-h)^{-1}\cap C \)]_{2d}
\\
	       &= [c(F-T)  \(c(T) s_*(Z(s),M_f)\)]_{2d}
\end{align*}
where we used the notation $h = c_1(\O_{\P(F \oplus \O)}(-1))$ and
\begin{equation}\label{Segreclass}
       s_*(Z(s),M_f) \eqdef \pi_* (1-h)^{-1} C
\end{equation}
for the total homological Segre class of the normal cone (this
terminology will be justified in a minute).
But $c(F-T)= c(\Ind Ds)^{-1}$
and since $T = TM_f$,
\begin{equation}\label{homchernclass}
       c_*(Z(s))\eqdef c(T) s_*(Z(s),M_f)
\end{equation}
is exactly the analogue of the
homological chern classes of \cite[example 4.2.6]{Fulton}.

We show that $c_*(Z(s))$ does not depend on the choice of $F$.	Again it
suffices to treat the case that $F'\subset F$.	We use primes whenever
an object is associated to $F'$.  The independence follows directly
from a formula for the Segre classes which expresses  how it behaves under
the extension $M'_f \subset M_f$ in terms of the  normal bundle
$F/F'$ of $M'_f \subset M_f$.
\begin{equation}\label{coneext}
       s_*(Z(s),M_f) = c(F/F')^{-1} s_*(Z(s),M'_f).
\end{equation}
Assuming \eqref{coneext}, we see that
\begin{align*}
       c_*(Z(s)) &= c(T)s_*(Z(s),M_f)
\\
		 &= c(T) c(F/F')^{-1}s_*(Z(s),M'_f) =
c(T')s_*(Z(s),M'_f).
\end{align*}
In particular we can take the limit over $K$.

Formula \eqref{coneext} is
well known for integrable complex manifolds \cite[example
4.1.5]{Fulton}, and we will follow the proof closely.
There are two terms in the class $C$ occurring in the definition
\eqref{Segreclass} of the Segre class, which we treat separately.

\nc\sfprimenul{\overline{(s'_f,0)}}
Note that there is a regular section $\sigma$
of $F/F'(1)$ on $\P(F\oplus\O)|_{M_f}$ cutting out
$\P(F' \oplus \O)|_{M_f}$. Therefore
\begin{align*}
       [\P(F'\oplus \O)|_{M'_f}]
	 &= e(F/F', s_f \bmod F') \cap [\P(F'\oplus \O)|_{M_f}]
\\
       &= e(F/F',s_f \bmod F') \cap e(F/F'(1),\sigma)
		       \cap [\P(F\oplus \O)|_{M_f}].
\end{align*}
Since on $\P(F'\oplus \O)|_{M_f}$ there is an exact sequence
$$
       0 \to Q'\to Q \to F/ F'\to 0,
$$
we have $e(Q',\sfprimenul)\cup e(F/F', s_f \bmod F') = e(Q,\sfnul)$. Then
the above implies that
\begin{align*}
	\Z(Q',\sfprimenul)
	 &= e(Q',\sfprimenul) \cap [\P(F'\oplus \O)|_{M'_f}]
\\
	&=e(Q,\sfnul )\cap e(F/F'(1),\sigma) \cap [\P(F\oplus \O)|_{M_f}]
\\
	&= e(F/F'(1),\sigma)\cap \Z(Q,\sfnul)
\end{align*}

As for the other term,
on $\~M_f$ there is a smooth section in $\O(-1)$ given
by $(s_f,0)$ which is an isomorphism $\O \iso \O(-1)$ on
$\~M_f -\E$. It follows that
$$
       [\~M'_f - \E] = e(F/F',s_f \bmod F') \cap [\~M_f - \E]
		      = e(F/F'(1),\sigma) \cap [\~M_f - \E].
$$
Then we have the equality
$$
       [\~M'_f] = e(F/F'(1),\sigma) \cap [\~M_f].
$$
because both left and right hand side are cycles supported on
$\~M'_f -\E \union \P(F'\oplus \O)|_{Z(s)}$ restricting to $[\M'_f
-\E]$.

For the computation of the Segre class we can forget
about the support given by $\sigma$ and use
$$
       e(F/F'(1)) = c_\Top(F/F'(1))
		   = \sum c_{\Top-j}(F/F')h^j.
$$
Thus we finally  get the expression
\begin{align*}
       s_*(Z(s), M'_f) &=
		\pi_*\(\sum h^{i+j} c_{\Top-j}(F'/F)
		    \cap (\Z(Q,\sfnul) - [\~M_f])\)
\\
		 &= c(F'/F) s_*(Z(s),M_f)
\end{align*}
which we set out to prove.

\bgroup
\def\N{\eN}
It remains to prove the claim~\eqref{dimclaim}.
We first turn to the case that $Z(s)$ is smooth
but possibly of the wrong
dimension.  This condition implies that $\Im Ds|_T \subset F$
has constant rank over $Z(s)$ because $\ker Ds|_{T}= \ker Ds = TZ(s)$.
Then $\Im Ds|_T$ is just the normal bundle $\N$ of $Z(s)$ in $M_f$.
Now let us identify the limit set $(s_f: 1/\lambda)(M_f)$ when $\lambda \to
\infty$. If we have a smooth path $\gamma$ with $\gamma(0) = x_0 \in
Z(s)$, then
we see that
$
       \lim_{t\to 0} (s_f:0)(\gamma(t)) = (Ds_f(\ddt|_0\gamma):0).
$
Therefore $\~M_f$ is just the blowup $\^M_f$ of $Z(s)$ in $M_f$.
This makes sense even though $M_f$ is only an almost complex manifold
since the normal bundle $\N$ has a complex structure.
The blow up is obtained
abstractly by identifying a tubular neighborhood $N_\epsilon$
of $Z(s)$ with the normal bundle, and replacing $N_\epsilon$ with
$I = \{(l, x) \in \P\N \times N_\epsilon \mid l \ni x\}$. It is  an
almost complex manifold, so certainly carries a fundamental class
$[\~M_f]$. It is also clear that $\E_f = \P\N$ is a submanifold
of real codimension $2$, and certainly satisfies the
claim~\eqref{dimclaim}.

Let $\O(\E_f)$ be the  smooth complex line bundle on the blow-up
$\^M_f$ defined by
the exceptional divisor $\E_f$, and let
$z \in A^0(\O(E))$ be a section cutting out $\E_f= \P\N$ with the
proper orientation i.e. $\Z(\O(\E_f),z) = [\E_f]$.
On $\^M_f$ the pulled back section is of the form $s_f = z \^s_f$ with
$\^s$ nowhere vanishing. Therefore the limit set of
$(s_f:1/\lambda)(\^M_f)$ in $\P(F \oplus \O)|_{\^M_f}$ as
$\lambda \to \infty$ is just $(\^s:0)(\^M_f) \union D$ where
$D \subset \P(F \oplus\O)|_{\E_f}$ is the $\P^1$ bundle
joining the zero section
$(0:1)|_{\E_f}$ and the section $(\^s_f:0)$.
Then down on $M_f$ the limit set of $(s_f: 1/\lambda)(M_f)$ is just
$\~M_f \union C\E_f$, where $C\E_f$ is cone bundle over $Z(s)$ joining $\E_f
\subset \~M_f$ and the zero section.

Now $C\E_f$ represents the homology class $C$. Thus
$$
       s_*(Z(s),M_f) = \pi_* (1-h)^{-1}C\E_f = \pi_* (1-h)^{-1} \E_f
		      = \pi_* (1-h)^{-1} \P\N = s(\N) \cap [Z(s)]
$$
Therefore if $Z(s)$ is smooth we find the expected formula
$$
       c_*(Z(s)) = c(TM_f)s(\N) \cap [Z(s)] = c(TZ(s)) \cap [Z(s)].
$$
Note that in deriving this formula we have not really used the
holomorphicity of $s$. It was sufficient that $M$ has an almost complex
structure and that $Ds$ is complex linear. Replacing manifolds by
stratified spaces the proof carries over essentially verbatim  if $Z(s)$
 is a local complete intersection since this condition implies that
$Ds|_{T}$ has constant rank, and that we have a well defined normal
bundle.
\egroup

In proving the claim \eqref{dimclaim} in the general case we  use
holomorphicity	more strongly.	We first blow up $Z(s)^{\mathrm{red}}$
in $M$ to get a new infinite dimensional  analytic space $\^M$. That
this is possible follows from the local analysis of the normal cone
in~\cite[\S III.1]{PT}.

Locally on $M$,   the exceptional divisor  $\E \subset \^M$ can be
described as follows.	Locally on $M$ we have an exact sequence of
holomorphic bundles
$$
       0 \to F^{\hol} \to E \to \~E^{\hol} \to 0,
$$
such that $TM|_{Z(s)} \surj\to \~E|_{Z(s)}$ is surjective, i.e. locally
$F^{\hol}$ can take the role of $F$.   Further, locally we can split the
sequence since $F^{\hol}$ has finite rank. Let the holomorphic subbundle
 $\~{\~E} \subset E$ be a lift of $\~E^{\hol}$. We write  $s =
s_f^{\hol} \oplus \~{\~s}$  corresponding to the decomposition $E =
F^{\hol} \oplus \~{\~ E}$.  Then locally $\E \iso \E^{\hol}_f
\times_{Z(s)} \P\~{\~E}$, where $\E^{\hol}$ is the exceptional divisor
of the blow up of $Z(s)$ in $M_f^{\hol}$, and where $M_f^{\hol}$ is the
integrable finite dimensional  complex manifold $Z(\~s^{\hol})$.
Moreover $\E_f^{\hol}$ is  naturally embedded in  $\P(F^{\hol} \oplus
\O)|_{Z(s)}  \subset \P(E \oplus \O)|_{Z(s)}$. If we are a little more
careful and choose $\~{\~E}$ such that	$\P\~{\~E}|_{Z(s)} \subset \E$
then  $\E = {\mathrm{Join}}(\E_f^{\hol}, \P\~{\~E}|_{Z(s)}) \subset \P
E|_{Z(s)}$.

Let $z \in H^0(\O(\E))$ be a section vanishing exactly along $\E$.
On $\^M$ we can decompose the section as $s = z^n \^s$.
Therefore, just as
in the previous finite dimensional case,
$(s: 1/\lambda)(\^M) \to \P(E \oplus \O)$ degenerates to $(\^s:
0)(\^M)\union n D$ where $D$ is the $\P^1$ bundle over $\E$ joining the
zero section $(0:1)|_{\E}$ and $(\^s_f:0)|_{\E}$.
Down on $M$, this means that $(s: 1/\lambda)(M) \subset \P(E \oplus \O)$
degenerates to $\~M \union C\E$ where $\~M \subset \P E$
is isomorphic to $\^M$ with $\~M \cap \P(E \oplus \O)|_{Z(s)} \iso \E$,
and $C\E$ is the cone bundle over $Z(s)$ joining the zero section and
$\E$.

Now we finally come to our claim~\eqref{dimclaim}.  The set
$\E_f = \P(F \oplus \O) \cap \E$.
At the very beginning we chose $F$
such that $F \supset F^{\hol}$.  Locally we define $\~{\~F} = F \cap
\~{\~E}$, then locally $F = F^{\hol} \oplus \~{\~F}$ and locally
$\E_f = {\mathrm {Join}}(\E^{\hol}_f , \P\~{\~F}|_{Z(s)})$. Thus $\E_f$ is
a stratified space of real dimension $2d +2N-2$, and we are done.
\end{pf}

\begin{Remark}\label{Zhat}
In the complex case we have obviously defined a class containing more
information about the section. Let
$$
       \widehat\Z(s) = c(\Ind(Ds))^{-1}c_*(Z(s)).
$$
\end{Remark}

\section{Seiberg Witten classes}

We will collect a few facts about Seiberg Witten basic classes in a
formulation suitable for arbitrary K\"ahler surfaces.
In the usual formulation, these classes
are the support of a certain function on the set of  $\Spin^c$-structures.
However in the presence of 2-torsion, $\Spin^c$-structures cause
endless confusion
which is why I have chosen to base my exposition on
SC-structures \cite{Karrer}. This notion catches the essence of
$\Spin^c$-structures, the existence of spinors.
It is  well suited to the K\"ahler case
and is equivalent to that of a $\Spin^c$-structure in dimension 4.
For more details see \cite{Karrer}. 

Let $X$ be a closed oriented manifold of dimension $2n$.
Choose a Riemannian metric $g$
with Levi-Civita connection $\nabla^g$,
and  Clifford  algebra bundle
$C(X,g) = C(T^\dual X,g)$.
There is a natural isomorphism of bundles
$c:\wedge ^* T^\dual X \to C(X,g)$ given
by anti-symmetrisation. It  induces a connection
and metric on $C(X,g)$ also denoted $\nabla^g$ and $g$.

An {\sl SC-structure}  is a smooth complex vector bundle $W$ of rank
$2^n$ together with  an algebra bundle isomorphism $\rho: C(X,g) \to
\sheafEnd(W)$.	In other words an SC structure is a bundle with the
irreducible Clifford algebra representation $\Delta$ in every fibre. A
section $\phi \in A^0(W)$ is called a (smooth) spinor.
An SC-structure exists if and only if $w_2(X)$ can be
lifted to the integers \cite[\S 3.4]{Karrer}. 
Existence will be clear in the case of K\"ahler surfaces.

SC-structures admit an invariant  hermitian metric i.e. one such that
Clifford multiplication by 1-forms is skew hermitian (sh). The chirality
operator $\Gamma = (\sqrt{-1})^n c(\Vol_g)$ has square $1$, and is
hermitian. Thus  $\Gamma$ has  an orthogonal  eigenbundle decomposition
$W = W^+ \oplus W^-$ with eigenvalue $\pm 1$, the positive and negative
spinors of the SC-structure. A one form $\omega \in A^1(X)$ defines  an
skew hermitian map $c(\omega): W^\pm \to W^\mp$  which is an isomorphism
away from the zero set	of $\omega$.

In this paragraph we assume $\dim(X)= 4$.
Then  $T^\dual_X \iso \sheafHom(W^+,W^-)^{\mathrm{sh}}$.
Let $L_W = \det W^+$. Then $L_W \iso \det W^-$, by the isomorphism
induced from
Clifford multiplication by a generic $1$-form,
which is an isomorphism outside codimension 4.
Thus $W$ is a $\Spin^c(4)$-bundle if we identify
$$
       \Spin^c(4) =
	\{(U_1,U_2) \in U(2) \times U(2) \mid \det(U_1) = \det(U_2)\}.
$$
We recover the usual definition $\Spin^c(4) = \Spin(4)\times_{\Z/2/Z}
U(1)$ from the isomorphism $\Spin(4) = SU(2) \times SU(2)$.
In any case  by chasing around the cohomology
sequences of the diagram
$$
\longmap
\begin{matrix}
0 &\map{}&\Z/2\Z&\map{}&\Spin^c(4)&\map{}&\SO(4)\times U(1)&\map{}1
\\
  &	 &\vequal{}&   &\upmap{}	  &	 &\upmap{}
\\
0 &\map{}&\Z/2\Z&\map{}&\Spin(4)  &\map{}&\SO(4)
&\map{}1
\end{matrix}
$$
we see that  $L_W + w_2(X) \equiv 0 \pmod 2$, and that this is the only
obstruction to lifting the $SO(4)\times U(1)$ bundle to $\Spin^c(4)$.
If $H^2(X,\Z)$
has no 2-torsion, the line bundle $L\equiv w_2(X)$  determines such a
lift completely,  and it is common to speak of the $\Spin^c$-structure $L$.

An {\sl SC-Clifford module} $(S,\<,>,\nabla)$,
is an SC-structure with a non-degenerate  invariant
hermitian metric $\<,>$ and
a unitary Clifford connection $\nabla$ i.e.  a unitary connection such
that for all vector fields $X$, spinors $\phi \in A^0(S)$, and $1$-forms
$\omega$ we have
$$
       \nabla_X (\omega \cdot \phi) = (\nabla^g_X \omega)\cdot \phi +
\omega \cdot \nabla_X \phi.
$$
The {\sl Dirac operator}  $\delbar$ of a Clifford module is the
composition
$$
       A^0(W) \map{\nabla} A^1(W) \map{\cdot} A^0(W).
$$
It is an elliptic self adjoint	first order differential operator, and
it maps positive spinors to negative ones and vice versa (i.e.
$\delbar:A^0(W^\pm) \to A^0(W^\mp)$).
Since $\rho$ is parallel, $\nabla$ respects the decomposition
$W = W^+ \oplus W^-$. Thus $\nabla$
induces a connection on $L_W$ with curvature $F$.

Much of the usefulness of SC-structures is a consequence of
the following easy lemma.

\begin{Lemma}
The set of isomorphism classes $\SC$ of SC-structures is an $H^2(X,\Z)$
torsor i.e. if $\SC \ne \emptyset$ and we fix an SC-structure $W_0$,
then for every SC-structures $W_1$, there exits a unique line bundle $\L$
such that $W_1 = W_0\tensor\L$.  Every SC-structure $S$  admits a
Clifford module structure $(W, \<,>, \nabla)$. If  we fix  one
SC-Clifford module $(W_0, \<,>_0, \nabla_0)$, there is a unique triple
$(\L,h,d)$  of a  smooth line bundle $\L$, with hermitian metric $h$ and
unitary connection $d$, such that
\begin{equation}\label{SCrepr}
       (W,\<,>,\nabla) \iso (W_0, \<,>_0, \nabla_0) \tensor (\L, h,d).
\end{equation}
\end{Lemma}

\begin{pf}
Clearly if $W_0$ is an SC structure, so is $W_0\tensor \L$ for every line
bundle $\L$. Conversely,
the bundle of Clifford linear homomorphisms
$\L(W_0,W) = \sheafHom_C(W_0,W)$ has rank~1, and the natural map $W_0
\tensor \L(W_0, W) \to W$ is an isomorphism.

For existence of a Clifford module structure see
\cite[prop. 4.2.1, 4.5.1]{Karrer}.
It will be clear for K\"ahler surfaces. It
follows directly from the definition of a Clifford module  that the
natural connection and metric on $\sheafHom(W_0,W)$ leaves $\L(W_0,W)$
invariant. Hence there is an induced metric and connection $(h,d)$ on
$\L(W_0,W)$, which has property \eqref{SCrepr}.  Conversely
if $(W,\<,>,\nabla)$ is defined by equation \eqref{SCrepr}, then
$$
       (\L,h,d) =
	  \sheafHom_C\((W_0,\<,>_0,\nabla_0)\,,\,
		(W_0,\<,>_0,\nabla_0)\tensor(\L,h,d)\)
$$
which proves uniqueness.
\end{pf}

If a base SC-structure is chosen, the line bundle $\L$ will be called
the twisting line bundle.


There is a natural gauge group $\G^\C$ acting on a Clifford module, the
group of all smooth invertible Clifford linear endomorphisms. $\G^\C$
can be canonically identified with $C^\infty(X, \C^*)$. In the
representation~\eqref{SCrepr}, $\G^\C = C^\infty(X, \C^*)$  acts in the
usual way on the set of metrics and unitary connections on the twisting
line bundle $\L$. Since every hermitian metric on a line bundle is gauge
equivalent, so is every Clifford invariant metric on a Clifford module.
Thus, up to gauge we can fix an invariant  metric and we are left with a
residual gauge group  $\G = C^\infty(X, \U(1))$.

The set of Clifford connections $\A$ on a fixed hermitian SC structure
$(W,\<,>)$ (i.e. Clifford module structures) is an affine space
modeled on $\sqrt{-1} A^1_\R(X)$.
Using the representation~\eqref{SCrepr} and harmonic representatives,
one shows that the set of connections mod gauge is
$$
       \B = \A / \G \iso \sqrt{-1} A^1_\R(X) / d\log C^\infty(X,\U(1))
	      \iso  H^1_{DR}(X)/H^1(X,\Z) \oplus \ker d^*
$$
We set $\Pee^* = \A \times A^0(W^+)^* / \G$. It is a
$\C\P^\infty \times \R^+$ bundle over $\B$. Thus $\Pee^*$ has the
homotopy type of $(S^1)^{b_1(X)} \times \C\P^\infty$.

There is an alternative description of $\B$ and $\Pee^*$ that will be
useful. Let $\A^\C$ be the set of all Clifford connections, and $\Herm$
the set of all hermitian metrics on $\L$. Let
$$
       \Amod  = \{(\nabla,<,>),\  \nabla
\txt{is} <,>\hbox{-unitary}\} \subset \A\times\Herm
$$
be the set of Clifford module structures.
Fix a metric $<,>_0$ and a $<,>_0$-unitary connection $\nabla_0$. The
representation	$\nabla = \nabla_0 + a$, models $\A^\C$ on $A^1_\C(X)$,
and the representation $<,> = e^f < , >_0$ models $\Herm$ on $A^0_\R(X)$.
A pair $(\nabla,<,> ) \in \Amod$ if and only if  $a+ \bar a = df$.
In particular $a$ is determined by $f$ and its imaginary  part,
so $\Amod$ is modeled on $A^0_\R(X)\times A^1_\R(X)$.

Now the diagonal action of $\G^\C$ on $\A^\C \times \Herm$
leaves $\Amod$ invariant. Our alternative description of $\B$ and
$\Pee^*$ is
\begin{equation}\label{alternative}
       \Pee^* = \Amod \times A^0(W^+)^*/ \G^\C \to \B = \Amod/ \G^\C
\end{equation}

Finally, to do	 decent gauge theory we have to complete to  Banach
spaces and -manifolds. Seiberg Witten theory works fine with  an $L^p_1$
completion of $\A$, $\A^\C$, and $A^0(W^+)$  and an $L^p_2$ completion
of $\G$, $\G^\C$   and $\Herm$ if $p > \dim X$. In this range $L^p_1
\inj\to C^0$, and therefore the  two possible $L_p$ descriptions of
$\Pee^*$ and $\B$ coincide. On the other hand,	 the Sobolev range does
not seem optimal: with more care and work one can probably use all
$p$-completions with $2 - \dim(X) / p > 0$. We will suppress completions
from the notation, explicitly mentioning completions if necessary.

{}From now on we assume $\dim X = 4$.
Fix an SC structure $W$ and choose an invariant hermitian metric $\<,>$.
Choose a Riemannian metric $g$ and a real  2
form $\epsilon$, which are {\em admissible} in the following sense:
$L_W$ admits no connection with  $F^+ = -2\pi \sqrt{-1}\epsilon^+$,
where as usual $+$ means taking the self dual part.
Admissible metrics and forms exist if $b_+ \ge 1$, since the condition is
certainly satisfied if
$c_1(L_W) \not\in \epsilon^{\mathrm{harm}} + H^-_g$ where $H^-_g$
is the space of $g$-anti-self-dual closed forms, and ``harm'' means
projection to the harmonic part.
Note that no use of Sard-Smale is made to define admissibility.
Actually for most of our purposes it would be  enough to
let $\epsilon$ be a closed (hence harmonic) self-dual form.

By a transversality argument \cite{Donaldson:intersectionform},
the admissible (metrics,forms) form  a connected set if $b_+ \ge 2$.
We say that a metric $g$ is admissible if
$(g,0)$ is.
Even if $b_+ = 1$, all metrics are admissible when $L_W^2 \ge 0$,
and $L_W$ is not torsion.

In dimension $4$,  the	anti-symmetrisation map gives an isomorphism
$c:\Lambda ^+ \iso \End_0^{sh}(W^+)$
between the real self-dual
forms and the traceless skew hermitian
endomorphisms of $W^+$. This special phenomenon allows us (or rather
Seiberg and Witten) to write down the  monopole equations \cite{Witten}
\begin{align}
\label{SW1}
       \delbar \phi &= 0  \qquad \phi \in A^0(W^+)
\\
\label{SW2}
       c(F^+)  &= 2\pi \phi\<\phi,-> - \pi |\phi|^2
		   -2\pi\sqrt{-1} c(\eps^+).
\end{align}
Let $\M = \M(W,g,\epsilon) \subset \Pee^*$
be the space of solutions modulo gauge.

%
%
%
As a technical remark, note that we use the conventions of \cite{BGV},
and that in their conventions  the Weitzenb\"ock (Lichnerowitz)
formula restricted to  $W^+$ reads
$$
       \delbar^2 = \nabla^* \nabla + r /4 + c(F^+/2)
$$
(\cite[th. 3.52]{BGV} and the observation that the
twisting curvature of an SC
structure is $1/\rank(W^+)$ times the curvature on $\det(W^+)$.)
The sign difference in the $c(F^+)$ term in \cite[lemma 2]{KM:Thom}
explains the relative change of
sign with respect to \cite[formula $(*)$]{KM:Thom}
in the Seiberg Witten equations. It is
chosen in such a way that the Weitzenb\"ock formula gives $C^0$ control
on the harmonic positive spinor $\phi$.

A basic property of the monopole equation noted by Witten,
which follows from the Weitzenb\"ock formula \cite[lemma 2]{KM:Thom} or a
variational description \cite[Section 3]{Witten}, is the following

\begin{Proposition}\label{Bochner}
The monopole equations have no solution with $\phi \ne 0$ if the
metric has positive scalar curvature.
\end{Proposition}

Alternatively we can define $\M$ as the zero of a Fredholm section in an
infinite dimensional vector bundle. Let
$\W^\pm = (\A\times A^0(W^+)^* \times_\G A^0(W^\pm) \to \Pee^*$.
Then $\M$ is the zero of the section in $\W^- \oplus A^+(X)$
given by the monopole equations~\eqref{SW1}, and \eqref{SW2}.

To see that it is actually a Fredholm section
we  linearise the equations, assuming
that  $(\nabla,\phi)$ is a solution, and $(\nabla + \eps a, \phi + \eps
\psi)$ with $a \in \sqrt{-1}A_\R^1(X)$ and $\psi \in A^0(W^+)$ is
a solution up to order 1 in $\eps$. We get (c.f \cite[eq.2.4]{Witten})
\begin{gather*}
	\delbar \psi + a\cdot \phi = 0
\\
	c^{-1}(2\pi(\phi\<\psi,-> + \psi\<\phi,-> - \Re\<\phi,\psi>)
-d^+a = 0.
\end{gather*}
The tangent space of the $\G$-orbit of $(\nabla,\phi)$ is
$\{(a,\psi) = (- d u, u \phi), \ u \in \sqrt{-1}A^0_\R(X)\}$.
Thus the Zariski tangent space of $\M$ in $(\nabla, \phi)$ is the first
cohomology of the Fredholm complex
$$
       \sqrt{-1} A_\R^0(X) \to \sqrt{-1}A_\R^1(X) \oplus A^0(W^+) \to
\sqrt{-1} A^+_\R(X)
\oplus
       A^0(W^-),
$$
where the maps are given by the left hand side of the linearised equations.
The virtual
dimension is given by Atiyah Singer index formula and is
\begin{equation}\label{vdim}
       d(W) = \vdim_\R(\M) = \quart(L_W^2 - (2 e(X) + 3 \sign(X))),
\end{equation}
where $e(X)$ is the topological Euler characteristic, and $\sign(X)$ the
signature \cite[eq. 2.5]{Witten}.

The crucial property that makes Seiberg Witten theory so much easier
than Donaldson theory is

\begin{Proposition}\cite[Corollary 3]{KM:Thom},\cite[\S 3]{Witten}
The moduli space $\M$ is compact.
For fixed $c >0$ there are only finitely  many SC-structure $W$ with
$d(\M(W)) \ge -c$  and $\M(W,g,\epsilon) \ne  \emptyset$.
\end{Proposition}

Note that for generic pairs $(g,\epsilon)$, moduli spaces of negative
virtual dimension are empty, but I do not see an a priori reason why
moduli spaces of arbitrary negative virtual dimension should not exist
for special pairs. Likewise for generic pairs the moduli space is
smooth of dimension $d(W)$ \cite{KM:Thom}. However we have no need
for this fact.

The index bundle $\Ind(Ds)$ of the deformation complex can be deformed
by compact operators (over a compact space !) into the sum of the index
of the signature complex and the index of the complex dirac operator.
Thus the determinant line bundle $\det(\Ind(Ds)$ of the index is
naturally oriented by choosing an orientation for $\det H^1(X,\R)^\dual
\tensor H^+(X,\R)$. We will in fact assume that an orientation for both
$H^+$ and $H^1$ is chosen. Suppose further that the pair $(g,\epsilon)$
is admissable (i.e. $\M((W,g,\epsilon) \subset \Pee^*$), then
proposition \ref{locEuler} in the previous section gives us  a
homology class $\MM \in  H_{d(W)}(\Pee^*)$,
i.e. a homology class
of the proper virtual dimension even if
$\M$ is not smooth, not reduced and not of the proper dimension (note that
in our case the moduli space $\M = Z(s)$ is compact, and  homology with
closed support is just ordinary homology).
In case $\M$  is smooth and has the proper dimension it is just the
fundamental class. The class $\MM$ depends only on the
connected component of $(g,\epsilon)$ in the space of admissable pairs,
by the homotopy property  of the localised Euler class proposition
\ref{locEuler}.\ref{homotopy}.
In particular $\M$ is independent of the admissable pair if $b_+ \ge 2$.

If $b_+ = 1$  the choice of an orientation of $H^+$ is the choice of a
connected component in
$\{\omega^2 >0\} \subset H^2(X,\R)$. It will be called the forward timelike
cone.  For every metric $g$ let $\omega_g$ be the unique self dual form
in the forward timelike cone with $\int \omega^2 = 1$. For a pair
$(g,\epsilon)$	and an SC-structure $W$ define the {\em discriminant}
\begin{equation}\label{discriminant}
    \Delta_W(g,\epsilon) = \int (c_1(L_W) - \epsilon) \omega_g
\end{equation}
A pair $(g,\epsilon)$ is admissable if	the
discriminant $\Delta_W(g,\epsilon) \ne 0$, because it means precisely
that $c_1(L_W) \notin \epsilon^{\harm} + H^-$. Clearly the discriminant
depends only on the period $(\omega_g, {\epsilon^+}^{\mathrm harm})$.

\begin{Lemma}
If  $b_+ = 1$ a pair $(g,\epsilon)$ is admissable if and only if
the discriminant $\Delta_W(g,\epsilon) \ne 0$.
There are exactly two connected components of admissable
pairs labeled by the sign of the discriminant.
\end{Lemma}

\begin{pf}
Suppose two pairs $(g_i,\epsilon_i)$,  $i=0,1$,
have discriminants $\Delta_i$  of equal sign.
Connect them by a  path $(g_t, \epsilon_t)$ in the space
of all pairs. Let $(\omega_t, \epsilon^{+,\harm}_t)$ be the corresponding
path of periods.
Then the discriminant
$$
       \Delta_t = \int (c_1(L_W) - \epsilon^{+,\harm}_t)\omega_t
$$
is continuous in $t$ but  may change sign.
However if we modify the path by setting
$$
      \epsilon'_t = \epsilon_t +
		  (\Delta_t - (1-t)\Delta_0- t\Delta_1)\omega_t
$$
then using
$\Delta_W(g,\epsilon + \delta) = \Delta_W(g,\epsilon) -  \int
\delta\wedge\omega_g$ and $\int \omega^2 = 1$ we see that
$$
 \Delta'_t= \Delta_W(g_t,\epsilon'_t) = (1-t)\Delta_0 + t\Delta_1.
$$
In particular $\Delta'_t$ does not change sign, so that
$(g_t,\epsilon'_t)$ is a path of admissable pairs.

Conversely if $c_1(L_W) \in \epsilon^\harm + H^-$, then any connection
$\nabla$ with induced Chern form $\epsilon^\harm$ determines a
``reducible''  solution $(\nabla,0) \in \Pee- \Pee^*$ of the monopole
equations.
\end{pf}

\begin{Definition}
If $b_+ \ge 2$, the {\em SW-multiplicity} is the map
\begin{align*}
       n:\SC &\to \Lambda^*H^1(X,\Z)[t] \iso
		  H_*(\Pee^*,\Z)
\\
       W &\maps\to \MM(W,g,\epsilon)
\end{align*}
where $(g,\epsilon)$ is any $W$-admissable pair.
If $b_+ = 1$ the {\em SW-multiplicities} $n_+$ and $n_-$ are defined
similarly but with pairs $(g_\pm,\epsilon_\pm)$ having positive
respectively negative discriminant.
\end{Definition}

It should be remarked that the SW-multiplicity (ies)  depend(s)
implicitly on the orientation of $H^+$ and $H^1$. For $b_+ > 1$ this is
only a matter of sign,	but for $b_+ = 1$ the orientation of $H^+$
determines in addition which invariant is $n_+$ and  which is $n_-$.

All known examples with $b_+ \ge 2$ have
non trivial multiplicities only when the virtual dimension
$d(W) = 0$. However
for surfaces with $p_g =0$ it is easy to give examples
with one of $n_\pm$ is non trivial for $d(W) >0$ we will in fact use such
an invariant. If $b_1 \ne 0$,
the $H^1$ part of the multiplicity becomes essential.

\begin{Remark}
Since $H_i(\Pee^*) = 0$ for $i <0$,
a moduli space of negative virtual dimension
never defines a nontrivial class.
Thus if for a class $L \in H^2(X,\Z)$ there exists an SC-structure $W$
with $L = c_1(L_W)$ and the multiplicity $n(W) \ne 0$ (respectively one
of $n_{\pm}(W) \ne 0$ then
$L^2 \ge 3 e(X) + 2\sign(X)$ (c.f. equation~\eqref{vdim}).
\end{Remark}

\begin{Remark}\label{specialchamber}
In the case $b_+ = 1$ we can alternatively consider the multiplicity as
depending in addition on a chamber structure in
$$
       \Gamma = \{(\omega,\epsilon) \in H^2(X,\R)^2
		\mid \omega^2 = 1,\ \omega_0 > 0\}
$$
where a chamber is defined by walls which are in turn defined by
all classes $L \equiv w_2(X)$ through equation~\eqref{discriminant}. This
is particularly useful when we consider structures  with $L_W^2 \ge
0$, $L_W$ is not  torsion. Then all pairs $(g,0)$ are admissable and
have discriminant of equal sign, because the forward timelike cone is
strictly on one side of the hyperplane $L_W^\perp \subset H^2(X,\R)$.
Thus for this subset we have a preferred chamber.
\end{Remark}

We will say that $L\in H^2(X,\Z)$ with $L \equiv w_2(X)$
has non trivial multiplicity if there is an SC-structure $W$ such that $L
= c_1(L_W)$ and $W$ has non trivial multiplicity. If $b_+ =1$ we will
further qualify which multiplicity is non trivial (i.e. $n_+$ or $n_-$)
or which chamber is chosen.  We will simply write $n(L) \ne 0$ or
$n_+(L) \ne 0$ etc.

A final and important piece of general theory is the following blow-up
formula \cite{Stern:talk},\cite[\S 8]{FS:rational}.
We will give a proof valid for K\"ahler surfaces in
section~\ref{computations}.

\begin{Theorem}\label{blowup4}
Let $X$ be a closed oriented 4-manifold with $b_+ \ge 1$.  An
SC-structure  $\~W$ on $X\#\Pbar^2$ can be decomposed as $\~W = W \#
W_k^{\Pbar^2}$, with determinant lines $L_{\~W} = L_W + (2k+1)E$.  If
the multiplicity $n_{(\pm)}(\~W) \ne 0$ then $d(\~W) = d(W) - k(k+1) \ge
0$, and the multiplicity $n_{(\pm)}(W) \ne 0$.	 Moreover if
$L_{W_{\Pbar^2}} = \pm E$ (i.e. $E \cdot L_{\~W} = \pm 1$)  then
$n_{(\pm)} (\~W) = n_{(\pm)}(W)$  under the identification $H^1(X,\Z)
\iso H^1(\~X,\Z)$.
\end{Theorem}

Here, $ n_{(\pm)} = n$ if $b_+ > 1$, and
if $b_+ = 1$, it is understood that we compare say  $n_+(W \#
W_k^{\Pbar^2})$ with  $n_+(W)$.

\section{Seiberg Witten classes of K\"ahler surfaces}

{}From now on, $(X, \Phi)$ denotes  a K\"ahler surface.
Then $X$ has a natural base SC-structure
$$
       W_0 = \Lambda^{0,*} X
$$
with Clifford multiplication given by
$$
       c(\omega^{10} + \omega^{01}) =
	    \sqrt2\(- i(\omega^{10}) +\eps(\omega^{01}) \),
$$
where $i$ is contraction and $\eps$ is exterior
multiplication. The metric and connection induced by the K\"ahler
structure on $\Lambda^{0*} X$ define a Clifford module structure on
$W_0$.
For an arbitrary SC structure $W$ = $W(\L)$ the spinor bundles
are of form
$$
       W^+ = (\Lambda^{00} \oplus \Lambda^{02})\tensor \L,
			      \qquad W^- = \Lambda^{01}(\L).
$$
and $L_W = \det(W^+) = -K \tensor \L^2$ (c.f. lemma~\ref{SCrepr}).
We call $\L$ the twisting line bundle.

We now turn to the monopole equations  (see also \cite[Section
4]{Witten}). In the decomposition of $W^+$, a positive spinor will be
written $\phi = (\alpha,\beta)$. The Dirac equation is then
\cite[Propos. 3.67]{BGV}.
$$
       \delbar \phi = \sqrt2 (\dbar\alpha + \dbar^* \beta)= 0.
$$


\nc\dzdz{dz_1\wedge dz_2}
\nc\dzbardzbar{d\bar z_1 \wedge d \bar z_2}
\def\dzdzbar#1{dz_#1\wedge d\bar z_#1}

Since $X$ is K\"ahler, we can locally choose holomorphic geodesic
coordinates $(z_1,z_2)$. A basis of the self dual forms is then the
K\"ahler form
$\Phi = \numfrac {\sqrt{-1}}2(\dzdzbar1 + \dzdzbar2)$,
$\dzdz$ and $\dzbardzbar$. Let $h$ be an hermitian metric on $\L$.
Choose a unit generator $e$ for $\L$,
then an orthonormal  basis for $W^+ $ is $e$ and
$\half e\dzbardzbar$.

Using the definition of Clifford Multiplication we compute:
\begin{align*}
       c(\Phi)e &= \numfrac{\sqrt{-1}}2
	  (-i(dz_1) \eps(d \bar z_1) +	\eps(d \bar z_1)i(d z_1)
	   -i(dz_2) \eps(d \bar z_2) +	\eps(d \bar z_2)i(d z_2))e
\\
		&= - 2\sqrt{-1} e.
\end{align*}
In exactly the same way we compute $c(\Phi)$, $\half e \dzbardzbar$,
and the action of $c(\dzdz)$  and $c(\dzbardzbar)$ on $e$ and
$\half e \dzbardzbar$. The result in matrix form is given by
$$
       c(\Phi) = \begin{pmatrix}
		      - 2\sqrt{-1} & 0
\\
			    0	  & 2\sqrt{-1}
		   \end{pmatrix}
\quad
       c(\dzdz) = \begin{pmatrix}
		0 & -4
\\
		0  & 0
		   \end{pmatrix}
\quad
       c(\dzbardzbar) = \begin{pmatrix}
		   0 & 0
\\
		   4 & 0
		   \end{pmatrix}.
\hskip 0pt minus 1 fil
$$

\def\bet#1#2{\beta_{\dot #1\dot#2}}
On the other hand, writing $\alpha = \alpha_e e$, and $\beta = \half
\bet12 e \dzbardzbar$,
$$
       (\alpha + \beta) \< \alpha +\beta,-> =
       \begin{pmatrix}
       |\alpha_e|^2   & \alpha_e\bar\bet12
\\
       \bar\alpha_e\bet12 & |\bet12|^2
       \end{pmatrix}.
$$
Thus if we define $\alpha^* = h(\alpha,-)$, $\beta^*= h(\beta,-)$ and
take the trace free part,
we get	the healthy global expression
$$
(2\pi(\alpha + \beta)\< \alpha + \beta, ->)_0=
 -2\pi\sqrt{-1} c\(\half (|\beta|_h^2 - |\alpha|_h^2) \Phi
	    + \sqrt{-1}(- \alpha \beta^* + \beta \alpha^*)) \)
$$

Plug all this in the monopole equations~\eqref{SW1},\eqref{SW2}.
Writing $c_1(F) =\numfrac{-1}{2 \pi i} F$, and using that $\Lambda \Phi
=2$ the monopole equation  for a K\"ahler metric and perturbation
$\epsilon = \lambda \Phi$ can be rewritten to
\begin{align}
       &\dbar \alpha + \dbar^* \beta = 0 \label{cSW1}
\\
       & F^{02}  =  2\pi \beta \alpha^*  \label{cSW2}
\\
       &F^{20}	=  -2\pi \alpha \beta^*   \label{cSW3}
\\
       &\Lambda c_1(F)^{11} = (|\beta|^2 - |\alpha|^2) + 2\lambda.
					   \label{cSW4}
\end{align}
Note  that $F$ is the curvature on $L_W$, but that   these are equations
for a unitary connection $d =\dee + \dbar$ on $\L$ and sections $\alpha
\in A^{00}(\L)$,  and $\beta \in A^{02}(\L)$ through the identity $F =
-F(K) + 2 F(\L,d)$. Here $F(K)$ is the	curvature of the canonical line
bundle i.e. minus the Ricci form.

In terms of the twisting bundle the virtual (real) dimension of the
moduli space reads
\begin{equation}\label{cvdim}
       d(\L) = d(\Lambda^{0*}(\L)) = \quart(L^2 - K^2)
	     = \L\cdot (\L- K).
\end{equation}
A more precise description is given by

\begin{Proposition} \label{Kahlermonopoles}
{\sloppy
A necessary condition for the existence of solutions
to the mono\-pole equations~\eqref{cSW1} to \eqref{cSW4}, is that
$(\L,\dbar)$ is a holomorphic line bundle, and that
}
\begin{align}
       -\deg_\Phi(K) &\le  \deg_\Phi(L) < \int(\lambda \Phi^2),
\txt{or}  \label{case0}
\\
	 \int \lambda \Phi^2 &< \deg_\Phi(L)  \le \deg_\Phi(K),
\txt{or}    \label{case2}
\\
       \int \lambda \Phi^2 &= \deg_\Phi(L)
	   \label{singcase}
\end{align}
In particular $L_W = -K \tensor \L^2$ has a natural holomorphic structure.
In case \eqref{case0}  the moduli space $\M = \M(\L, \Phi, \lambda)$ of
solutions can be
identified as a real analytic space with the moduli space of pairs of a
holomorphic structures $\dbar$ on $\L$, and a divisor $\alpha \in
|(\L,\dbar)|$, in particular the Zariski tangent space in
$(\dbar,\alpha)$ is canonically identified with
$H^0(\L|_{Z(\alpha)})$.
In case \eqref{case2} the moduli space $\M$ of solutions
can be identified with the moduli
space of pairs of a holomorphic structure $\dbar$ on $\L$, and an
element $\beta \in \P H^2(\L) = |K\tensor \L^\dual|^\dual$, in
particular the Zariski tangent space at $(\dbar,\beta)$ is isomorphic to
$\overline{H^0(K\tensor \L|_{Z(\bar\beta)})}$.
In case \eqref{singcase} the ``moduli space'' $\M \subset \Pee - \Pee^*$
(i.e. $\alpha = \beta = 0$)
can be identified with the space of holomorphic structures $\dbar$ on $\L$.
\end{Proposition}

\begin{pf}
Combining~\eqref{cSW1} and ~\eqref{cSW2} yields
\begin{equation}\label{positivity}
       \dbar\dbar^* \beta = -\dbar^2 \alpha
			   = - F^{02}\alpha
			   = -2\pi |\alpha|^2 \beta.
\end{equation}
Integrating both sides against $\<\beta,->$, immediately gives that
$\alpha\beta= 0$ and $\dbar\beta = \dbar\alpha = 0$. Thus $F^{02} =
F^{20} =0$, Since $F^{02} = 2 F^{02}(\L,d)$, $\dbar$ is a
holomorphic structure on $\L$, and either $0 \ne \alpha \in H^0(\L)$ and
$\beta= 0$ or $0 \ne \beta \in H^2(\L)$ and $\alpha = 0$,  or $\alpha =
\beta = 0$. Note
that if for example $\alpha \ne 0$, then $\beta = 0$ is cut out
transversely by equation~\eqref{positivity}.
The last monopole equation~\eqref{cSW4} gives the condition
$$
       \deg(L) = -\deg(K) +  2\deg(\L) = \half\int\Lambda c_1(F)\Phi^2
	     = \half\int ( |\beta|^2 - |\alpha|^2  + 2\lambda) \Phi^2
$$
which fixes the global $L_2$ norm of $\alpha$ and $\beta$, and determines
whether $\alpha \ne 0$ or $\beta \ne 0$ or $\alpha = \beta = 0$.

Finally we deal with equation~\eqref{cSW4}.  If $\alpha \ne 0$ and
$\beta =0$ then we are dealing essentially with the abelian vortex
equation studied by Steve Bradlow \cite[\S 4]{Bradlow}
Oscar Garcia-Prada
and earlier in a different guise by Kazdan Warner \cite{KazdanWarner}.
See also \cite{Bradlow:nonabelian} and \cite{OkonekTeleman:coupledSW}.
I thank Steve Bradlow for pointing out	that almost all of the work had
already been done by him and Oscar Garcia-Prada.
To identify the moduli space as a real analytic space we just
jazz up Bradlow's results a bit. This is necessary because we have to
to understand how
the moduli space is cut out in order to apply the localised Euler
class machinery in the next section.

It is slightly more convenient to use our alternative
description~\eqref{alternative} of $\Pee^*$, and solve for a pair
$(d_L, h)$ where  $h = e^f h_0$ is a hermitian	metric on $L$ and $d_L=
\dee + \dbar = d_0 + a$ is $h$ unitary, and mod out the  full gauge
group $\G^\C$ of all complex nowhere vanishing functions.  To be
precise we take $d_L$ in $L^p_1$, and $\G^\C$ and $f$ in $L^p_2$  with
$p > 4$. The sections $\alpha$ and $\beta$, being disguised spinors, are
as before in $L^p_1$.

For an $h$-unitary connection we have , $\dee h(s,t) = h(\dee s, t) +
h(s, \dbar t)$ for all sections $s,t \in A^0(\L)$. Thus $d_L$ is
determined by $\dbar$ and $h$, or equivalently, $a^{01}$ and $f$.
Expressed in $a^{01}$ and $f$, equation~\eqref{cSW4} becomes
\begin{equation}\label{masterf}
   \laplace f = 2\pi (|\beta|^2_{h_0} - |\alpha|_{h_0}^2)e ^f  -
2\sqrt{-1}\Lambda(\dee_0 a^{01} - \dbar_0\bar{a^{01}}) +   \mu
\end{equation}
where
$\mu = 2\pi(2\lambda + (\Lambda c_1(F(K)) - 2 \Lambda c_1(\L,\nabla_0)$
(compare \cite[lemma 4.1]{Bradlow}).

If $\beta$ is small in $L^p_1$ hence in $C^0$,
we can solve for $f$ in equation~\eqref{masterf} with
the solution depending real analytically on $(a^{01},\alpha)$ by
the analytical lemma~\ref{fsoln}.
Moreover, variation of ~\eqref{masterf} with respect
to $f$ when $\beta = 0$ gives
\begin{equation}\label{deltamasterf}
    \delta\hbox{``eqn \eqref{masterf}''} =
	       (\laplace  +  2\pi|\alpha|^2e^f)\delta f.
\end{equation}
Thus, equation~\eqref{masterf} cuts out this solution
transversely. More invariantly, if $\beta$ is small, there is a unique
metric $h(\dbar,\alpha,\beta) = h_0e^{f(\dbar-\dbar_0,\alpha,\beta)}$
solving the last  monopole equation~\eqref{cSW4}.

\begin{Lemma} \label{fsoln}
Let $X$ be a compact Riemannian manifold,
and $\dim(X) < p < \infty$ a Sobolev weight.
Then for every real non negative function
$0 \le w_0 \in L^p$,  with $\int w_0 > 0$
and real function $\mu_0 \in L^p$, with $\int \mu_0 > 0$, there exists
a neighborhood $U_{(w_0,\mu_0)} \subset L^p \times L^p$  such that
for all $(w,\mu)\in U_{(w_0,\mu_0)}$ the equation
\begin{equation}\label{rawfeq}
       \laplace f = - w e^f + \mu
\end{equation}
has a unique $L^p_2$ solution depending analytically on $w$ and $\mu$.
The solution is smooth if $w$ and $\mu$ are smooth.
\end{Lemma}

\begin{pf}
As in \cite[lemma 4]{Bradlow} make the substitution
$f = \~f - g$ where $g$ is the unique solution of $\laplace g  =
\int \mu - \mu $ to reduce to the case where $\mu$ is constant.
Then apply \cite[theorem 10.5(a)]{KazdanWarner}
to solve the equation for $w_0,\mu_0$ (note that Kazdan Warners Laplacian
is negative definite and that the proof
works fine with $w \in L^p$ instead of
$C^\infty$).
Since at a solution $f_0$ for $(w_0,\mu_0)$ we have
$$
       \delta \hbox{``eqn \eqref{rawfeq}''} =
	      (\laplace + w_0 e^{f_0})\delta f
$$
and $(\laplace + w_0 e^{f_0})$ is invertible, we conclude with the
implicit function theorem that there continues to exist  a solution  for
$(w,\mu)$ in a small neighborhood of $(w_0,\mu_0)$, and that this
solution depends real analytically on $(w,\mu)$. Regularity follows from
standard bootstrapping techniques. Uniqueness follows  from the weak
maximum principle (\cite[theorem 8.1]{GilbargTrudinger},  c.f.
\cite[remark 10.12]{KazdanWarner}).
\end{pf}

In geometric terms, this has the following consequence. let $\A^{01}$ be
(the $L^p_1$-completion) of the space of $\dbar$-operators on $\L$
modeled on $A^{01}(X)$ through $\dbar = \dbar_0 + a^{01}$. The complex
gauge group $\G^\C$ acts naturally by conjugation. Let
$$
       \Pee^{01*} = \A^{01} \times (A^{00}(\L)\oplus A^{02}(\L))^* /\G^\C
$$
Clearly there is a projection
$\Pee^* \to \Pee^{01*}$ forgetting $h$. What we have done is
showing that there is section
\begin{align*}
       \Pee^{01*} & \to  \Pee^*
\\
       (\dbar,\alpha,\beta) &\to (\dbar,\alpha,\beta,h(\dbar,\alpha,\beta))
\end{align*}
in a neighborhood of $\beta = 0$,
whose image is cut out as a real analytic space by the last monopole
equation~\eqref{cSW4}.

So far we have not used the other equations. Suppose we are in
case~\eqref{case0}, i.e. where a solution corresponds to sections.
Then $\M$ is cut out by $\dbar^2 =0$, $\dbar \alpha
= 0$,  $\beta = 0$ and, by the preceding argument, $h =
h(\dbar,\alpha,\beta)$. Thus projection identifies $\M$
with
$$
       \MBN = \{ (\dbar,\alpha,\beta) \in \Pee^{01},
		    \ \dbar^2 = 0,\ \dbar\alpha =0, \beta=0\}
$$
For the Zariski tangent space it gives
\begin{align*}
       T_{(\nabla,\alpha,0,h)}\M &= T_{(\dbar,\alpha,0)}\MBN
\\
			    &= \Ker\left.
					 \begin{pmatrix}
				  \dbar  & \alpha \\
						 & \dbar
					  \end{pmatrix}
				     \right/ \Im
					  \begin{pmatrix}
					     \alpha \\
					     -\dbar
					   \end{pmatrix}
\\
				 &= {\Bbb H}^1( 0 \map{} \O
					 \map{\alpha} \L \map{} 0)
\\
			   &= H^0(\L|_{Z(\alpha)}).
\end{align*}
It is easy to check that the linearised versions of equations~\eqref{cSW1},
\eqref{cSW2}, \eqref{cSW3},
and~\eqref{masterf} give the same result (as it should).

\comment{
Substituting $(\nabla + \eps a,\alpha + \eps \xi,\eps  \eta, h e^{\eps
f})$  in \eqref{complexSW}, when  $(\nabla,\alpha,0,h) \in \M$ yields
\begin{align*}
       &\dbar^*\eta + \dbar \xi + a^{01} \alpha = 0
\\
       &\dbar a^{01} = \bar\alpha \eta
\\
       &\dee a^{10} = - \alpha \bar \eta
\\
       &\numfrac{\sqrt{-1}}{2\pi}\Lambda (\dbar a^{10} + \dee a^{01}) +
\numfrac1{4\pi} \laplace f = -|\alpha|^2_h f  - 2 \Re
h(\alpha,\xi).
\end{align*}
As above, the equations give $\eta = 0$ and we can just solve for $f$.
We mod out the	infinitesimal gauge transformations sending $u \in
A^0_\C(X)$ to
$(a,\xi,\eta,f) = (-du,u\alpha,0,2\Re u)$
The conclusion is the same.
}

Case \eqref{case2} is reduced to the previous case by Serre duality.
In case \eqref{singcase} the metric $h$ we look for is an
(almost) Hermite-Einstein metric.
\end{pf}


\begin{Corollary} \label{genus}
Let $X$ be K\"ahler surface.
and $L \equiv w_2(X)$ be a  class in $H^2(X,\Z)$ with $n(L) \ne 0$.
Then $L$ is of type $(1,1)$.
Moreover if $p_g >0$, then for all K\"ahler forms $\Phi$ on
$X$, the class	$L$ satisfies
$$
       \deg_\Phi(K_X) \ge \deg_\Phi(L) \ge -\deg_\Phi K_X
$$
If $p_g = 0$, and $n_-(L) \ne 0$ (resp. $n_+(L) \ne 0$), then
$$
       \deg_\Phi(L) \ge  -\deg_\Phi(K_X)\
		\txt{(resp.} \deg_\Phi(L)\le \deg_\Phi(K_X))
$$
\end{Corollary}

\begin{pf}
First we consider the case $p_g >0$. Under the conditions of the
corollary, there is an SC-structure $W$ with $L_W = L$ which admits at
least one  solution to the monopole equation for {\em every} admissable
pair $(g,\epsilon)$. In particular $W$ admits a solution for every
K\"ahler metric and $\epsilon = \lambda \Phi$. Thus $L = L_W$ is of type
$(1,1)$. Moreover the necessary condition for the existence of a
solution of section {\em or} cosection type (i.e. equation  \ref{case0}
{\em or} \ref{case2} in proposition~\ref{Kahlermonopoles}) gives
precisely the required inequality in the limit	$\lambda \to 0$.

If $p_g =0$, then $L$ is automatically of type $(1,1)$ and
say the condition $n_-(L) \ne 0$ means that there is an SC structure $W$
with $L_W = L$ such that for  any  K\"ahler metric,
 $W$ admits solutions  of section type (i.e. equation
\ref{case0}) if  $\lambda$ is sufficiently large.
This gives a lower bound but no upper bound on $\deg_\Phi(L)$.
\end{pf}

\begin{Remark}\label{easyp_g=0}
If $p_g =0$ and we restrict to perturbation $\epsilon =0$ (or small),
then  the same argument as in the $p_g >0$ case gives the stronger
degree inequality if  $L^2 \ge 0$, $L$ is not torsion,	since in this
case all metrics are admissable and  have discriminant of equal sign. In
particular on a Del Pezzo surface such classes do not exist.
\end{Remark}

\begin{Corollary} \label{Kisthere}
Let $X$ be a K\"ahler surface with base SC structure $W_0 =
\Lambda^{0*} X$.  Then $n(W_0) = 1$ if $p_g >0 $ and $n_-(W_0) = 1$ if
$p_g =0$, in particular $n(-K_X) \ne 0$ resp. $n_-(-K_X) \ne 0$.
Likewise, $n(W_0(K_X) = \pm 1$ if $p_g>0$ and $n_+(W_0(K_X) = \pm 1$,
in particular $n(K_X)\ne 0$ resp. $n_+(K_X) \ne 0$.  Moreover $W_0$ is
the only SC-structure $W$ with $L_W = -K_X$ mod torsion and non
trivial multiplicity $n$ respectively $n_-$. In particular if $L \in
H^2(X,\Z)$, such that $L = -K \in H^2(X,\Q)$ and $n(L) \ne 0$ resp.
$n_-(L) \ne 0$ then $L = -K \in H^2(X,\Z)$.
\end{Corollary}

\begin{pf}
We will prove the statement for $-K_X$.  Then we have to consider
SC-structures $W = \Lambda^{0*}(\L)$ with $c_1(\L)$ torsion. Choose a
K\"ahler metric and $\lambda \gg 0$. Then $\M(W) \iso  \MBN(\L)$ the
moduli space of line bundles with a section. But $\M^{BN}(\L)$	 is
just a reduced point if $\L$ is trivial, and empty if $c_1(\L)$ is  non
trivial torsion. Thus $W_0 = \Lambda^{0*} X$ is unique among the
SC-structures $W$ with $L_W = -K_X$ mod torsion  with $n(W) \ne 0$
(resp. $n_-(W) \ne 0$). In fact its multiplicity is $1$.
The case $+K_X$ can be dealt similarly with Serre duality.
Its multiplicity is $\pm 1$ because of the unpleasant orientation
switches.
\end{pf}

\begin{Corollary}\label{slickdivisor}
Let $D$ be an  effective divisor with $D\cdot(D-K) = 0$, $h^0(\O(D))
=1$,  $h^0(\O_D(D)) =0$, and $h^0(\L(D)) = 0$ for every line bundle $\L
\in \Pic^0(X)$. Then $n(-K_X + 2D) \ne 0$ if $p_g >0$  and $n_-(-K_X +
2D) \ne 0$ if $p_g =0$. Likewise, $n(K_X -2D) \ne 0$ if $p_g >0$ and
$n_+(K_X - 2D) \ne 0$ if $p_g = 0$.
\end{Corollary}

\begin{pf}
This corollary is proved just as the previous one, and reduces to it if
$D =0$. The conditions of the corollary  ensure  precisely  that
$\MBN(\O(D))$ consists of one smooth point and that
$\vdim(\Lambda^{0*}(D)) = 0$.
\end{pf}

We are finally in the  position to prove the main theorem~\ref{main}
and corollary~\ref{poscurv}.

Our first task is to define a set $\K$ of basic classes.

\begin{Definition}
If $b_+ \ge 2$ then the basic classes are defined by
$$
       \K = \{ K \in H^2(X,\Z) \mid n(K) \ne 0\}
$$
If $b_+ = 1$ then $\K = \K_- \union \K_+$ where
\begin{align*}
       \K_- =\{ K &\in H^2(X,\Z) \mid n_-(K) \ne 0, \txt{and}
	     \exists L \txt{with} n_-(L) \ne 0
\\
       &\txt{such that} n_-(L-m(K+L))\ne 0 \txt{for some $m \ge 1$} \}.
\end{align*}
The set $\K_+$ is defined similarly in terms of $n_+$. Here we are allowed
to take $m\ge 1$ rational as long as $m(K+L)$ is two divisible.
\end{Definition}

These basic classes
are rightfully {\em the} Seiberg-Witten basic classes when $b_+
\ge 2$, but for $b_+ =1$ the definition is geared towards the specific
application we have in mind.
We will show that $\K$ has all properties~\ref{*}.

It is clear that $\K$ is an oriented diffeomorphism invariant, and that
the basic classes are characteristic.  The pushforward property
\ref{*}.\ref{iii} follows immediately from the blow up formula
theorem~\ref{blowup4} or~\ref{blowup}. For K\"ahler surfaces
the classes are of type $(1,1)$ by
corollary~\ref{genus}.
The degree property~(\ref{*}.\ref{ii})
(for all surfaces minimal or not) follows also from
corollary~\ref{genus}. This is immediate for $p_g >0$.
If $p_g = 0$ assume that $K \in \K_+$ say, the case $K \in \K_-$ being
essentially the same.
Then the corollary gives the three inequalities
\begin{align}
       \deg K &\le  \deg K_X,
\\
	\deg L &\le  \deg K_X,
\\
	-m\deg K &\le  \deg K_X + (m-1) \deg L \le m\deg K_X.
\end{align}
If $p_g >0$ then $K_X \in \K$ by corollary~\ref{Kisthere}.
Thus it remains to check that $K_X \in \K$ if $p_g =0$.
In fact we will check that $-K_X \in \K$.

We have already seen in corollary~\ref{Kisthere} that $n_-(-K_X) \ne 0$.
Either directly from corollary~\ref{slickdivisor},
or using the invariance under the reflection in the
exceptional curves $E_1, \ldots, E_n$ we see that
$n_-(-K_X + 2\sum E_i) \ne 0$. Then
denoting
$$
       \L_m = m\Kmin + \sum E_i,
$$
we have to check that $n_-(-K_X + 2\L_m) \ne 0$.
We will distinguish  four  cases.


If $\kod(X) = 0$, then $\Kmin$ is torsion and
we can take $m=\rmmath{ord}(\Kmin)$, since $n_-(-K_X + 2\sum E_i) \ne 0$.

If $\kod(X) = 1$, then $\Xmin$ has a unique elliptic fibration	$\Xmin
\map{\pi} C$. By the canonical bundle formula,	$\Kmin =
\pi^*\L_C(\pi^*K_C + \sum (p_i - 1)F_i)$, where $\L_C$	is a line bundle
on $C$ of degree $\chi$. Since $p_g = 0$ and $\chi \ge 0$,  we have $0
\le g\le q\le 1$, and we distinguish further between $g=0$ and $g =1$.

If $g= 0$, then $c_1(\pi^*\L_C(K_C)) =(\chi-2) F$,  where $F$ is a
general fibre, and there are at least $3-\chi$ multiple fibers because
$\Kmin >0$. Now the class  $\Kmin + \sum_{i=1}^{2-\chi} F_i =
\sum_{j=3-\chi}^n (p_j -1) F_j$ is of the form	$m\Kmin$ with rational
$m >1$. Again by corollary~\ref{slickdivisor}, we have
$$
       n_-(-K_X + 2\L_m) = n_-(-K_X + 2(\sum_{j=3-\chi}^n (p_j-1)F_j +
\sum E_i)) \ne 0
$$

If $g=1$, then $\chi = 0$, and $K_C =0$.
In this case we can  take $m = 1$ since $c_1(\L_C) =0\in H^2(X,\Z)$
and by corollary~\ref{slickdivisor}
$$
   n_-(-K_X +2\L_1) = n_-(-K_X + 2(\sum (p_i-1)F_i + \sum E_i)) \ne 0.
$$

The most instructive case is when $X$ is of general type. Then the
irregularity $q=0$ since $p_g = 0$ and $\chi(\O_X) >0$.
Take  $m=2$,
then $\MBN(\L_2) = |2\Kmin + \sum E_i|$. By
formula~\eqref{gtpluri} (or directly by Ramanujan vanishing)
$$
       \dim_\C\MBN(\L_2) = P_2 - 1 = \Kmin^2  = \half \vdim_\R(W_2).
$$
Thus the moduli space is again smooth of the proper dimension and we
conclude that $n_-(-K_X + 2\L_2) \ne 0$.   In fact
$n_-(\Lambda^{0*}(\L_2)) = t^{\Kmin^2}$ since the $\O(1)$ on $\Pee^*$
corresponds to the $\O(1)$ on $\MBN$. This is because both measure the
weight of the action of the constant gauge transformations on the
spinors respectively sections.

It now follows	from lemma~\ref{inequality}  that if $\kod(X) \ge 0$,
all SW-structure  have a moduli space of
virtual dimension $d=0$,  and up to torsion, the basic classes are of
type
\begin{equation}
       K =\lambda \Kmin + \sum \pm E_i \bmod
       \hbox{Torsion}, \qquad |\lambda| \le 1.
\end{equation}
Moreover by proposition~\ref{Kcharacterisation},  $\Kmin$ is invariant
up to sign and torsion and   every $(-1)$-sphere is represented by a
$(-1)$-curve up to sign and torsion.

We first get rid of  torsion in the $(-1)$-curve conjecture i.e. theorem
{}~\ref{main} part~\ref{mainb}.
Let $e$ be a $(-1)$-sphere, giving a connected sum decomposition $X =
X'\# \Pbar^2$.
As we have used before, there is a
diffeomorphism $R_e =\id \# \C$-conjugation  representing the reflection
in $e$.

I claim that for any SC-structure $W$ on a 4-manifold
$$
       R_e^*(W) = W \tensor\O((c_1(L_W),e)e),
$$
where $\O(e)$ is the
line bundle corresponding to the Poincar\'e dual of $e$.
In fact if we write $R_e^* W = W \tensor \L$, then
$\L = \sheafHom_C(W, R_e^*W)$ (c.f. the proof of \ref{SCrepr}).
Now we can just identify $W$ and $R_e^* W$ on $X'$, i.e. $\L$ is
trivialised on $X'$. Thus
$$
       c_1(\L) \in \Im H^2(X, X- X', \Z) \iso
	 H^2(\Pbar^2) \subset H^2(X,\Z).
$$
Write  $\L = \O(a e)$ for some integer $a$. Since
$$
       L_W + 2 a e = L_{R_e^* W} = R_e^*L_W = L + 2(e,L_W) e
$$
the claim is proved.

Going back to the K\"ahler case, we can assume that  $e$ is homologous
to a $(-1)$-curve $E$ up to torsion. Consider $W =
R_e^*R_E^*(\Lambda^{0*} X) = \Lambda^{0*}(E-e)$. By oriented
diffeomorphism invariance $n_{(-)}(W) \ne 0$  (in case $p_g =0$ we have
tacitly used the fact that $R_e^*R_E^*$ induces the identity on rational
cohomology so in particular does not change the orientation of $H^+$).
Moreover  $c_1(L_W) = -K_X$ up to torsion. By corollary~\ref{Kisthere},
we conclude  that $W = \Lambda^{0*} X$, so $e = E \in H^2(X,\Z)$.

Finally for the invariance $\pm \Kmin$, consider any basic class of the
form $L = \pm \Kmin + \sum \pm E_i$ up to torsion. After reflections in
the $(-1)$-curves, we get a class $L'$ equal to  $\pm K_X$ up to
torsion. By corollary~\ref{Kisthere}  $L' = \pm K_X \in H^2(X,\Z)$. Now
for any basis $E'_1, \ldots E'_n$ of the lattice in $H^2(X,\Z)$ spanned
by the $(-1)$-spheres (e.g. the $(-1)$-curves) we have the identity
$$
       \pm \Kmin  = L' + \sum (E'_i,L') E'_i
		  = L + \sum ( E'_i, L) E'_i \in H^2(X,\Z).
$$
This finally proves theorem~\ref{main}.

\begin{Remark}
It is easy to give a
definition of basic classes for $b_+ = 1$
that satisfies all properties ~\ref{*}	except the
invariance under blow down (i.e. property~\ref{*}.\ref{iii}).
A class $K$ is then basic if there exists a metric $g$ such
for all $\delta >0$
there exists an admissable pair $(g,\epsilon)$ with
$\|\epsilon^{+,\harm}\| < \delta$ such that
$n(g,\epsilon,K) \ne 0$. The degree inequality for minimal surfaces
then follows from remark \ref{specialchamber}. But alas, if $K^2 <0$ one
can not avoid the possibility that a chamber on the blow up
realisable with small $\epsilon$
can only be realised for large $\epsilon$ on the blow down.
In my original treatment I used this  definition. I am grateful to
Robert Friedman whose
insistent questions about my definitions made me realise this mistake.
\end{Remark}

\begin{Remark}\label{whycastelnuovo}
An  easy application of the
techniques of the next section gives the following.
If $\L$ is a holomorphic line bundle  on a surface with $p_g = q = 0$ with
$h^0(\L) \ge \chi(\L)\ge 1$,
then $n_-(\Lambda^{0*}(\L)) = t^{\L(\L-K_X) \over 2}$.
If $p_g= q = 0$ and $\kod(X) \ge 0$
we can apply this to $\L_2= 2\Kmin + \sum E_i$. Then by the Castelnuovo
criterion  and the above  we conclude
$n_-(-K_X + 2L_2) \ne 0$. This gives an alternative way to prove
that $-K_X \in \K$  in	this case.
%
%
Conversely the degree inequality~\ref{*}.\ref{ii} cannot hold true for
rational and ruled surfaces for K\"ahler forms $\Phi$ such that
$\deg_\Phi(K_X) < 0$. Since in deriving the degree inequality we
did not use that $\kod(X) \ge 0$, we conclude that for $\kod(X) =
-\infty$ the set of the above defined basic classes $\K= \emptyset$.
In particular we see that the following proposition is a rather
direct analog of to the classical  Castelnuovo criterion.
\end{Remark}

\begin{Proposition}\label{Castelnuovo}
A K\"ahler surface
is rational if and only if $b_1= 0$, and $\K = \emptyset$.
\end{Proposition}

\begin{Remark}
After reading \cite{FM:SW} I realised the following. The blow up
formula~\ref{blowup4} can be generalised to connected sum decompositions
$X = X'\# N$ with $N$ negative definite and $H_1(N,\Z) = 0$. The latter
condition is automatic for K\"ahler surfaces of non negative Kodaira
dimension by a beautiful observation of Kotschick
(an unramified covering $\~N \to N$ of
degree $d$ gives an unramified covering $\~X = d X' \# \~N \to X'\# N$
which is an algebraic surface of non negative Kodaira dimension with a
connected sum decomposition with a factor with $b_+  >0$). Such smooth
negative definite manifolds $N$ have $H_2(N) = \directsum_{i=1}^n \Z
n_i$. SC structures $W_N$ on $N$ are determined by $L_N = \sum (2a_i +
1) n_i$. Thus the reflections $R_{n_i}$ in $n_i^{\perp}$,  act on the SC
structures on $N$. SC -structures on $X'\# N$ are of the form $W= W_{X'}
\# W_N$. Now the blow up formula is as if $N = n \Pbar^2$:   $W = W_{X'}
\# W_N$ is an SW-structure on $X'\# N$ if and only if $W_{X'}$ is a
SW-structure on $X'$ and $d(W) \ge 0$. In particular the Seiberg Witten
structures are invariant under the operation $R_{n_i}: W_{X'}\# W_N \to
W_{X'} \# R_{n_i} W_N$, and $\sheafHom_C(W, R_{n_i}W)$ has a
trivialisation over $X'$.  With these remarks the arguments for
$(-1)$-spheres carry over directly to prove that  for K\"ahler surfaces
$X$ with  $\kod(X)\ge 0$, with a connected sum decomposition $X = X'\#
N$, $H_2(N) \subset H_2(X)$ is spanned by $(-1)$-curves.
\end{Remark}

Stefan Bauer showed me how to use the  Seiberg Witten multiplicities and
the basic classes to  determine the multiplicities of the elliptic
surface. If the surface does not have finite cyclic fundamental group,
the multiplicities can be read off from the topology.	Thus we consider
a minimal elliptic surface $X_{pq}$ fibred over $\P^1$ with 2 multiple
fibers of multiplicity $p$ and $q$ We will assume  that  $p\le q$.

\begin{Corollary} (Bauer)
The multiplicities $p$ and $q$ are determined by the underlying oriented
differentiable manifold, unless $p_g = 0$, $p=1$ and $q$ arbitrary.  The
surfaces $X_{1q}$ are all rational and diffeomorphic.
\end{Corollary}

\begin{pf}
If the canonical class $K_X$ is not  torsion , we can
write $K_X$ in terms of the primitive vector $\kappa$ in
the ray spanned by $K_X$, normalised so that $\kappa \Phi >0$
$$
       K_X = (p_g - 1)F + (p-1)F_p + (q-1) F_q
	  = {(p_g + 1)pq - p -q \over \gcd(p,q)} \kappa
		\in H^2(X,\Z)/\mathord{\hbox{Torsion}}.
$$
Let  $d(p,q) = \((p_g+1)pq-p-q\)/\gcd(p,q)$  be the oriented
divisibility of $K_X$. If $K_X$ is torsion we simply set $d(p,q) = 0$.

The divisibility $d(p,q)<0$ if and only if $p_g=0$, $p=1$ and $q$ is
arbitrary. But this implies that $K_X$ is rational. We have already seen
that we can recognise rationality as Kodaira dimension $-\infty$ and
$b_1 =0$ (corollary~\ref{Kodaira} or proposition~\ref{Castelnuovo}).
Thus we can assume that $X_{pq}$ has non negative Kodaira dimension.
Then $\pm K_X \in \K$ are the basic  classes  with  the highest
divisibility (or torsion) and the oriented divisibility  $d(p,q) \ge 0$
is just the  unoriented divisibility of $\pm K_X$. The number
$\gcd(p,q)$ is also determined by the oriented manifold, being the order
of the fundamental group. Choose one of these classes, say $-K_X$.

First consider the case $p_g >0$.  Suppose that  $K= -K_X + 2 F_q \le
0$, (i.e. on the same side of $0$ as  $-K_X$), then it is the basic
class with second largest divisibility since $F_q$ is the smallest
effective vertical divisor, and  $n(-K_X + 2F_q)) \ne 0$ by
lemma~\ref{slickdivisor}  above. Thus if there exist basic classes other
then $\pm K_X$, we can reconstruct $p$ from $(2p/ \gcd(p,q))\kappa = K -
(-K_X)$. Since $d(p,q)$, $p_g$ and $\gcd(p,q)$ are known, this
determines $q$ as well. Obviously if we have chosen $+K_X$  the same
arguments works with  $K = K_X - 2F_q$,  there is nothing that prefers
$K_X$ over $-K_X$.

In the case $p_g =0$ we make a small modification. We choose an
orientation of $H^+$, which for a moment we assume is the standard one.
Consider  the classes $K \in H^2(X,\Z)$ mod torsion in the half ray
spanned by $0$ and  $ -K_X$ with unoriented divisibility at most
$d(p,q)$ (i.e. in between $0$ and $-K_X$)  such that $n_-(K) \ne 0$.
Note that $-K_X$ is just the basic class with largest divisibility in
$\K_-$.  Then if $K = -K_X + 2 F_q \le 0$ we can use exactly the same
argument as in the case $p_g >0$.
If we choose a different orientation of $H^+$, we replace
$-K_X$ by $+K_X$ but just as above the conclusion is the same.

If $\K	= \pm K_X$ or
for $p_g =0$ if $\{K \in [-K_X, 0] \mid n_-(K) \ne 0\} = -K_X $
then  $d(p,q)\gcd(p,q)< 2p$.
The few possibilities are listed in the following table
$$
\begin{array}{|l|c|c|c|l|}
\hline
\strut&(p,q) &\gcd(p,q) &d(p,q) &\text{Type}	 \\
\hline
p_g=0 &(2,2) &2 	&0	&\text{Enriques} \\
      &(2,3) &1 	&1	&		 \\
      &(2,4) &2 	&1	&		 \\
      &(2,5) &1 	&3	&		 \\
      &(3,3) &3 	&1	&		 \\
      &(3,4) &1 	&5	&		 \\
\hline
p_g=1 &(1,1) &1 	&0	&\text{K3}	 \\
      &(1,2) &1 	&1	&		 \\
\hline
\end{array}
$$
Clearly, in this case the pair $(p,q)$ is determined by the
oriented differentiable manifold as well.
\end{pf}

To prove that no surface with $\kod \ge 0$ admits a metric with positive
scalar curvature (corollary~\ref{poscurv}), first consider the case $p_g
>0$. Then the statement is clear, and one of Witten's basic
observations.  By proposition~\ref{Bochner}, for 4-manifolds with positive
scalar curvature  $n(K) = 0$ for all $K \in H^2(X,\Z)$, since
for our metric with positive scalar curvature $g$ and	 small
perturbations $\epsilon$, we have $\M(W,g,\epsilon) = \emptyset$ for all
SC-structures $W$. On the other hand we just showed that $n(-K_X) \ne 0$
using a K\"ahler metric.

The same argument works if $p_g = 0$ and $K_X^2 \ge 0$:
$n(-K_X,g,\epsilon)$ is independent of the metric $g$ and of $\epsilon$ as
long as $\epsilon$ is small, with the exception of the case $-K_X$
torsion in which case we have to choose  $\epsilon$ in the forward light
cone. But we can do better.

For the general case $p_g =0$, we choose a perturbation  $\epsilon =
\lambda \Phi$ with $0<\lambda \ll 1$ say.  Now suppose that the metric
with positive scalar curvature $g$ has period $\omega_g = \omega_{\min}
+ \sum \eta_i E_i$ where $\omega_{\min} $ is the projection to the
cohomology of minimal model.  Then since $\omega_g$ is in the interior
of the forward light cone,  and $\Kmin$ is in the closure of the forward
light cone, $\omega\cdot \Kmin = \omega_{\min}\cdot \Kmin\ge 0$ with
equality iff $\Kmin$ is torsion. Then for {\em some} choice of signs in
$-\Kmin - \sum \pm E_i$ we have
$$
  \omega_g\cdot (-\Kmin - \sum \pm E_i)\le 0 < \lambda \int\omega_g \Phi
$$
Thus for {\em some} choice of signs we compute $n_-$ (rather than $n_+$)
with our metric of positive scalar curvature and small perturbation.
Hence	 $n_-(-\Kmin-\sum \pm E_i) = 0$. On the other hand   $n_-(-\Kmin
-\sum \pm E_i) = n_-(-K_X) \ne 0$,  a contradiction just like before.

\section{Some computations of Seiberg-Witten multiplicities}%
\label{computations}

In this section we will go beyond determining potential basic classes and
compute the Seiberg Witten multiplicity of elliptic surfaces. We also
prove an algebraic version of the  blow up formula.
It is here that our excess intersection formulas pay off.
We first show how to go over to a fully complex point of view.
Then we use the special geometry of elliptic surfaces to compute the
multiplicities and finally we prove a blow up formula.

{}From now on we identify an SC-structure with the corresponding
twisting line bundle $\L$.
We will consider the solutions of the monopole equations of section
type, i.e. corresponding to equation~\eqref{case0}.

We have already seen  that the variation of the last monopole
equation~\eqref{cSW4} with respect to the hermitian metric is $h$ is
given by  $(\laplace + |\alpha|_h^2)h^{-1}\delta h$ (c.f.
equation~\ref{masterf}).  Therefore the  solutions to the fourth
monopole equation ~\eqref{cSW4} is a smooth submanifold of $\Pee^*$ in a
neighborhood of  the moduli space $\M(\L)$. In the proof of
proposition~\ref{Kahlermonopoles} we have seen that we can identify this
submanifold with the   ``vortex locus'' $\{h = h(\dbar,\alpha,\beta)\}$
i.e. the image of the section $\Pee^{01*} \to \Pee^*$. The vortex locus
is well defined in a neighborhood of the moduli space $\M(\L)$ only, but
this will not affect our arguments, as the construction of the localised
Euler class in section \ref{top} depends only on a neighborhood of
$\M(\L)$. Since the vortex locus is given by a function, we can identify
it with its domain $\Pee^{01*}$ which carries a natural complex
structure.

By property \ref{stability} of	proposition~\ref{locEuler} we are
allowed to compute the localised Euler class   $\MM(\L)$ of the moduli
space by considering $\M(\L)$ as a zero set of a section $S$ over  the
vortex locus cut out by the remaining equations, which define the same
ideal as
$$
       \dbar^2 = 0, \qquad \dbar \alpha =0, \quad \beta = 0
$$
i.e.  complex  equations ! Moreover the deformation complex of these
equations on
$\Pee^{01*}$ in a point $(\dbar,\alpha,0)$ is
$$
       A^{00}(X) \map{}
	A^{01}(X) \oplus A^{00}(\L) \oplus A^{02}(\L) \map{}
	A^{02}(X) \oplus A^{01}(\L)
$$
where the map is complex linear. We trivialise the determinant of the
index using the complex structure. This has brought us safely in
complex waters, and allows us to use proposition ~\ref{locChern} and in
particular formula~\ref{magic}.

{}From now on we identify $\M(\L)$ with $\MBN(\L)$.
Define the vector bundles
$$
       \sfA^{pq}(\L) =\(\A^{01} \times
     (A^{00}(\L \oplus A^{02}(\L))^*\)\times_{\G^\C} A^{pq}(\L)
$$
over $\Pee^{01*}$. Then $\MBN$ is given by a section $s$ in $E = A^{02}(X)
\oplus \sfA^{01}(\L)$, and the tangent space is given by
$$
       T\Pee^{01*} \iso \(A^{01}(X) \oplus \sfA^{00}(\L)
			       \oplus \sfA^{02}(\L)\)/ A^{00}(X).
$$
The  deformation complex can be considered as a map
$T\Pee^* \to  E$ and is exactly what we called $Ds$ in section~\ref{top}.

To identify the index $\Ind Ds$ we first make a compact
perturbation, keeping only the differential operator part of the
deformation complex. Then it  splits naturally in the $\dbar$ complex on
$X$ with index $\C^{\chi(\O_X)}$ and the index of complex
$$
       0 \map{} \sfA^{00}(\L) \map{\dbar} \sfA^{01}(\L)
				 \map{\dbar} \sfA^{02}(\L) \to 0
$$
where $\dbar$ is the universal $\dbar$ operator descended to $\Pee^{01*}$.

To rewrite this index
in holomorphic terms, consider the universal divisor
$$
       \Delta = \{ (\dbar,\alpha,x) \mid \alpha(x) = 0\}
$$
on the pull back of $X \times \MBN$.  Now if $\Omega^{pq}$ is the sheaf
of $C^{\infty}$ $(p,q)$-forms on $X$ considered as an $\O(X)$ module,
then I claim that on the pull back of $X \times \MBN$ to $\A^{01} \times
(A^{00}(\L \oplus A^{02}(\L))^*$,  there is a $\G^\C$ equivariant  exact
sequence
$$
       0 \map{} \O(\Delta)
	  \map{\dbar} p_1^*\Omega^{00}(\L)
	  \map{\dbar} p_1^*\Omega^{01}(\L)
	 \map{\dbar} p_1^*\Omega^{02}(\L) \to 0.
$$
In fact this only says that $(\dbar,\alpha, x) \to \alpha(x)$ is a
$\G^\C$ invariant section vanishing along $\Delta$ with multiplicity $1$
lying in the kernel of $\dbar$, which is obvious.  Now	descend this
whole complex to $X \times \MBN$ and take push forward to $\MBN$.  Then
we get an exact sequence of complexes
$$
       0 \map{} Rp_*\O(\Delta) \map{\dbar} \sfA^{00}(\L)
	 \map{\dbar} A^{01}(\L) \map{\dbar} \sfA^{02}(\L) \map{} 0
$$
where we are considering $Rp_*\O(\Delta)$) as a complex with zero
boundary operator and $\sfA^{pq}(\L)$ as a complex concentrated in
degree $0$. Thus for the index we find
\begin{equation}\label{index}
       \Ind(Ds)  = \Ind\(Rp_* \O(\Delta)\)  + \C^{\chi}
\end{equation}

A more precise description of $\MBN$  depends on the surface.
Here we will do the  case of elliptic surfaces. The author has succeeded
in treating ruled surfaces in a similar way.

\begin{Proposition}\label{multiplic}
Let $X \map{\pi} C$ be a K\"ahlerian elliptic surface
over a curve $C$ of genus $g$, with
multiple fibers $F_1, \ldots F_r$ of multiplicity $p_1, \ldots p_r$
of holomorphic Euler characteristic $\chi$.
Consider the line bundle $\L = \O(\pi^*D + \sum a_i F_i)$ where
$D$ is a divisor on $C$ of degree $d$, and $0 \le a_i < p_i$.
Then the Seiberg Witten multiplicity is zero if $d < 0$, and if $d\ge 0$
it is given by
$$
       n_{(-)}(\Lambda^{0*}(\L)) = \begin{cases}
	     (-1)^d{\chi + 2g-2\choose d}  & \txt{if} \chi + g -2  \ge 0
\\
       \sum (-1)^j { 1-g-\chi + d-j \choose d-j}{g \choose j}
      & \txt{if} \chi + g - 2 < 0
      \end{cases}
$$
\end{Proposition}

Note that if the topological Euler characteristic $e >0$ (or equivalently
$\chi >0$) then $g =q = \half b_1(X)$ \cite[corollary II.2.4]{FM}, so in
this case  $\chi + g-2 = p_g  - 1$. Note that the second formula is just
1 if $p_g =q =0$ (i.e. $e>0$). This illustrates
remark~\ref{whycastelnuovo}.

If $p_g >0$ and $q = g = 0$, so in particular $e = 12 \chi > 0$, Witten
proves this formula by choosing a general $\omega \in H^0(K_X)$ and
using  the perturbation $\epsilon = \omega + \bar\omega$. He then argues
that the multiplicity $n(\L)$ is the number of ways we can decompose a
fixed canonical divisor $K_0$ as  $K_0 = D_+ + D_-$ with $D_+\in
|(\L,\dbar_0)|$, and $D_- \in |K\tensor(\L,\dbar_0)^\dual|$, where
$\dbar_0$ is the unique holomorphic structure that $\L$ admits \cite[eq.
(4.23) e.v.]{Witten}.

To be honest,  this is
what I read out of it. Note for example that his sheaf $R$ is just
$\L|_Z(\alpha)$, and that
$$
       h^0(R) = h^0(\L|_{Z(\alpha))} = \dim T_{(\dbar,\alpha)} = d
$$
(the last equality we will see in a minute).
Actually I think that the
computations below are the mathematical version of (I paraphrase)
``integrating over
the bosonic and fermionic collective coordinates in the path integral''
and ``computing the Euler class of the bundle of the cokernel of the
operator describing the linearised monopole equations over the moduli
space (the bundle of
antighost zero modes)'' \cite[above (4.11)]{Witten}. In fact with
hindsight, the latter seems  a dual description
of the localised Euler class in the case that the
cokernel has constant rank.

\begin{pf}
We choose a K\"ahler metric and $\lambda$ such that
$\deg_\Phi(\L^{\tensor 2}(-K)) < \lambda \Vol(X)$.  This means that if
$\L$ has non zero multiplicity, it must carry a holomorphic structure
with a section. In case $p_g =0$ it also means we are looking at $n_-$.
But $(\L,\dbar)$ has a section if and only if $D$ is an effective
divisor on $C$. In fact a family of vertical line bundles with a section
gives a family	of effective divisors on $C$ by pushforward of the line
bundle, and conversely a family of effective divisors on $C$ gives a
family of vertical line bundles with a section	by pull back and
multiplication with a fixed section in $\O(B) = \O(\sum a_i F_i)$ ($B$
for base locus). Thus there is a natural isomorphism
$$
       \MBN \iso \MBN_C = C^d
$$
where $C^d$ is the $d^{\text{th}}$ symmetric power of $C$. The functorial
isomorphism comes with an isomorphism $\O(\Delta_X) = \O(\pi^*\Delta_C +
B)$.

Next we use Grothendieck Riemann Roch (an alias of the family index
theorem). Let $q: C \times C^d	\to C^d$ be the projection map. Then the
projection $p: X \times \MBN \to \MBN$ can be factored as  $p = q\circ
\pi\times \id$.  Thus writing $\pi\times \id$ as $\pi$,
\begin{align*}
       \ch(Rp_*\O(\Delta))
       &= \ch\(Rq_* \(\O(\Delta_C) \tensor R\pi_*\O(B)\)\)
\\
       &= q_*\(ch(\O(\Delta_C) \ch R\pi_*\O(B) \td(C)\)
\\
       &= q_*\(ch(\O(\Delta_C))
	   \pi_*\(e^B(1 - K/2 + \chi(\O_X)(pt\times C^d)\)\)
\\
       &= \chi(\O_X) q_*\(ch(\O(\Delta_C))(pt\times C^d)\)
\\
       &= \ch(\O(1)^\chi),
\end{align*}
where we have abbreviated the holomorphic Euler characteristic by
$\chi$. If we denote by $x$ the chern class of $\O(1)$, then our
computation shows that
$$
       c_t(\Ind Ds) = (1 + tx)^\chi,
$$
at least over the rationals.

The chern classes of the tangent bundle of $C^d$ are computed in
\cite[eq. VII.5.4]{ACGH}. Denoting the pullback of the $\theta$ divisor
on $\Pic^d$ to
$C^d$ by $\theta$ the result is
$$
       c_t(T_{C^d}) = (1+ tx)^{d+1-g} e^{-t\theta / (1+ tx)}
$$
Combining  these two expressions, our multiplicity drops out
\begin{align*}
       n(\Lambda^{0*}(\L)) &= c(\Ind Ds)^{-1} c(T_C^d) \cap [C^d]
\\
	     &= [(1 + tx) ^{d +1-g - \chi} e^{-t\theta /1+tx}]_{t^d}
\end{align*}
With the following identity of formal power series
\cite[eq. VIII.3.1]{ACGH}
$$
       [(1 + xt)^a f( - t/(1+ xt))]_{t^b}  = [ (1  -  xt)^{b - a - 1}
f(-t)]_{t^b}.
$$
the expression	becomes
$$
       n(\Lambda^{0*}(\L)) = [(1-tx)^{\chi + g - 2} e^{-t\theta}]_{t^d}
=
\begin{cases}
(-1)^d\sum_{j=0}^d {\chi + g - 2 \choose d-j} {\theta^j x^{d-j} \over j!}
		   & \txt{if} \chi + g- 2 \ge 0
\\
\sum_{j=0}^d  (-1)^j {1-g-\chi + d-j \over d-j}{\theta^j x^{d-j} \over j!}
	     & \txt{if} \chi + g -2 < 0
\end{cases}
\hskip 0pt minus 1fil
$$
Now $\theta^jx^{d-j}\cap [C^d] = j! {g \choose j}$ \cite[below eq.
VIII.3.3]{ACGH}. The elementary identity $\sum_j {a \choose j}{b \choose
c-j} = {a + b \choose c}$ then gives the answer as stated.
\end{pf}

As a second application of the methods developed we give
a complex analytic version of the blow up formula.

\begin{Proposition}\label{blowup}
Let $(X,\Phi)$ be a K\"ahler surface, and $\L$ a line bundle on $X$.
Suppose that $\deg_\Phi(\L^{\tensor 2}(-K) < \lambda \Vol(X)$.	Let
$\sigma:\~X \to X$ be the blow up of $X$ in a point,
with K\"ahler form $\~\Phi$, and let
$\~\L = \L(a E)$ be a line bundle on $\~X$ with $a\ge 0$. Suppose
that the cohomology class of $\~\Phi$ is close to $\Phi$.
 Then  there is a natural identification $\M(\~\L) =
\M(\L)$, and
$$
       \MM(\~\L) = [(1+x)^{a(a-1)/2}\widehat \MM(\L)]_{\dim_\R=
\L\cdot(\L-K)- a(a-1)}.
$$
Here $\widehat\MM$ is the class defined in remark~\ref{Zhat}, and $x$
the class of the natural bundle $\O(1)$ over $\M$.
In particular if $a=0,1$ then $n(\~\L) = n(\L)$.
\end{Proposition}

Of course this proposition determines the multiplicity
$$
       n_{(-)}(\Lambda^{0*}(\L(aE))) = n_{(-)}(\Lambda^{0*}(\L(-aE)))).
$$
Since quite in general $n_+(\Lambda^{0*}(\L)) = \pm n_-(\Lambda^{0*}(K
\tensor \L^{\dual})$ it determines the corresponding relation for
$n_+$ up to sign which is really all we need here.

\begin{pf}
The conditions on the degree imply that a solution of the monopole
equations correspond to a holomorphic structure on $\L$ with a section.
Since $\~\Phi$ is close to $\Phi$ we have (by definition of close)
$\deg_{\~\Phi}(\~L) < \lambda\Vol(\~X)$, hence solutions on the blowup
also correspond to holomorphic structures on  $\~\L$ with  a section.

Now $aE$ is contained in the base locus of the sections. Thus similar to
what we did for elliptic surfaces, we get an identification of $\M(\L)$
with $\M(\~\L)$ by multiplication of the
section with a section in $\O(aE)$, and the universal divisor
on $\~X\times \M(\~\L)$ is $\~\Delta = \Delta + aE$.

Again, identify the chern class of the index of the
deformation complex with formula~\eqref{index}.
Let $\~p$ be the projection $\~X \times \M(\L) \to \M(\L)$, and
$p$ the projection $X\times \M(\L) \to \M(\L)$.
Then the total chern class of the index is
$$
   c(R\~p_*(\~\Delta)) = c\(Rp_*\(\O(\Delta)\tensor R\sigma_*\O(aE)\) \).
$$
By induction on $a$,  one shows that
$$
       R\sigma_* \O(aE) = \O -\O_{pt}^{a(a-1)/2}.
$$
Since $\O(\Delta|_{pt\times \M(\L)}) = \O(1)$ it gives
$$
       c(R\~p_*(\~\Delta)) = c(Rp_*\O(\Delta))/c(\O(1))^{a(a-1)/2}.
$$
Formula~\eqref{magic} gives us
$$
       \MM(\~L) =
	[(1+x)^{a(a-1)/2} \(c(Rp_*(\Delta))^{-1}c_*(\M(\L)\)]_{d(\~L)}
$$
Since the real virtual dimension of $\M(\~\L)$ is
$d(\~\L) = \L\cdot (\L-K) -a(a-1)$ and the term in brackets is exactly
$\widehat\MM(\L)$, we have proved the formula.
\end{pf}

\nopagebreak
\end{document}